\documentclass[11pt]{article}
\usepackage{amscd}
\usepackage{amssymb,amsfonts,amsmath,psfrag,amscd,cite}
\usepackage{bbm}
\usepackage{graphicx}

\newtheorem{theo}{Theorem}
\newtheorem{prop}[theo]{Proposition}
\newtheorem{coro}[theo]{Corollary}
\newtheorem{defi}[theo]{Definition}

\newtheorem{rema}[theo]{Remark}

\makeatletter \@addtoreset{equation}{section}
\@addtoreset{theo}{section} \makeatother

\begin{document}

\title{ Poisson Reduction of Controlled Hamiltonian System by
Controllability Distribution }
\author{Tudor S. Ratiu\\
Section de Math\'{e}matiques and Bernoulli Center,\\
\'{E}cole Polytechnique F\'{e}d\'{e}rale de Lausanne. CH-1015 Lausanne. Switzerland\\
E-mail: tudor.ratiu@epfl.ch\\\\
Hong Wang \thanks { Corresponding author Hong Wang (E-mail: hongwang@nankai.edu.cn).} \\
School of Mathematical Sciences and LPMC,\\
Nankai University, Tianjin 300071, P.R.China\\
E-mail: hongwang@nankai.edu.cn }
\date{ \emph{ Dedicated to Professor Manuel de Le$\acute{o}$n
on the Occasion of His 60th Birthday }\\
December 25, 2013} \maketitle
 
{\bf Abstract.} In this paper, we first study the Poisson reductions
of controlled Hamiltonian (CH) system and symmetric CH system by
controllability distributions. These reductions are the extension of
Poisson reductions by distribution for Poisson manifolds to that for
phase spaces of CH systems with external force and control. We give
Poisson reducible conditions of CH system by controllability
distribution, and prove that the Poisson reducible property for CH
systems leaves invariant under the CH-equivalence. Moreover, we
study the Poisson reduction of symmetric CH system by $G$-invariant
controllability distribution. Next, we consider the singular Poisson
reduction and SPR-CH-equivalence for CH system with symmetry, and
prove the singular Poisson reduction theorem of CH system. We also
study the relationship between Poisson reduction for singular
Poisson reducible CH systems by $G$-invariant controllability
distribution and that for associated reduced CH system by reduced
controllability distribution. At last, some examples are given to
state the theoretical results.\\



\tableofcontents
\section{Introduction}

Poisson reduction is a important topic in the study of Poisson
geometry for constructing new Poisson manifold and simplifying
Hamiltonian systems defined a Poisson manifold, and it is also a
powerful tool in the study of stability and bifurcation theory of
mechanical systems. See Abraham et al \cite{abma78,abmara88},
Libermann and Marle \cite{lima87}, Marsden et al \cite{mamiorpera07,
ma92, mawe74}, Marsden and Ratiu \cite{mara99}, Ortega and Ratiu
\cite{orra04}. Just as we have known that the Poisson reduction is
first a generalization of symplectic reduction method to Poisson
manifolds and to the singular context. In addition, one can also
study the Poisson reduction for Poisson manifolds by pseudo-groups
and distributions, since the special Poisson brackets can be induced
on the presheaf of Poisson algebras. There have been many results
and ways of reduction, such as, optimal point and orbit Poisson
reduction, regular Poisson reduction for Hamiltonian systems by
using optimal momentum map, are given in Ortega and Ratiu
\cite{orra04}, and the singular reductions for Hamiltonian system,
Dirac structure and nonlinear control system, are given in Sjamaar
and Lerman \cite{sjle91}, Jotz et al \cite{jorasn11},
$\acute{S}$niatycki \cite{sn06}, as well as the reduction for
Poisson manifolds by distributions, are given in Marsden and Ratiu
\cite{mara86}, Falceto and Zambon \cite{faza08}, and Jotz and Ratiu
\cite{jora09}, and
there is still much to be done in this topic.\\

On the other hand, it is well known that the mechanical control
system theory has formed an important subject in recent twenty
years. Its research gathers together some separate areas of research
such as mechanics, differential geometry and nonlinear control
theory, etc., and the emphasis of this research on geometry is
motivated by the aim of understanding the structure of equations of
motion of the system in a way that helps both analysis and design.
So, it is natural to study mechanical control systems by combining
with the analysis of dynamic systems and the geometric reduction
theory of Hamiltonian and Lagrangian systems. For examples, Birtea
et al in \cite{bipura04} and S$\acute{a}$nchez de Alvarez in
\cite{sa89}, studied the controllability of Poisson systems;
Nijmeijer and Van der Schaft in \cite{nivds90}, studied the
nonlinear dynamical control systems as well as the use of feedback
control to stabilize mechanical systems; Leonard and Marsden in
\cite{lema97} and Bloch and Leonard in \cite{blle02} gave the
underwater vehicle with internal rotors and rigid spacecraft with
internal rotor as the practical models of Hamiltonian systems with
control. In particular, we note that in Marsden et al
\cite{mawazh10}, the authors studied regular reduction theory of
controlled Hamiltonian systems with symplectic structure and
symmetry, as an extension of regular symplectic reduction theory of
Hamiltonian systems under regular controlled Hamiltonian equivalence
conditions, and Wang in \cite{wa12} generalized the work in
\cite{mawazh10} to study the singular reduction theory of regular
controlled Hamiltonian systems, and Wang and Zhang in \cite{wazh12}
generalized the work in \cite{mawazh10} to study optimal reduction
theory of controlled Hamiltonian (CH) systems with Poisson structure
and symmetry by using optimal momentum map and reduced Poisson
tensor (or reduced symplectic form). In addition, since the
Hamilton-Jacobi theory is developed based on the Hamiltonian picture
of dynamics, it is natural idea to extend the Hamilton-Jacobi theory
to the (regular) controlled Hamiltonian system and its a variety of
reduced systems, and it is also possible to describe the
relationship between the CH-equivalence for controlled Hamiltonian
systems and the solutions of corresponding Hamilton-Jacobi
equations. Wang in \cite{wa13d, wa13f, wa13e} studied this work and
applied to give explicitly the motion equations and Hamilton-Jacobi
equations of reduced spacecraft-rotor system and reduced underwater
vehicle-rotors system on a symplectic leaves by calculation in
detail, which show the effect on controls in regular symplectic
reduction (by stages) and Hamilton-Jacobi theory.\\

It is worthy of note that if there is no momentum map for our
considered systems, then the reduction procedures given in Marsden
et al \cite{mawazh10} and Wang and Zhang \cite{wazh12}, and hence in
Wang \cite{wa12,wa13d,wa13f,wa13e} can not work. One must look for a
new way. On the other hand, motivated by the work of Poisson
reductions by distribution for Poisson manifolds in Ortega and Ratiu
\cite{orra04}, we note that in above research work, the phase space
$T^*Q$ of the CH system is a Poisson manifold, and its control
subset $W\subset T^*Q$ is a fiber submanifold. If we assume that
$D\subset TT^*Q |_W$ is a controllability distribution of the CH
system, then it is a natural problem if we could study the Poisson
reduction for the CH system by controllability distribution. This is
our first goal of research in this paper. Next, Wang and Zhang in
\cite{wazh12} give the regular Poisson reduction of CH system, we
hope to develop this reduction in the singular context, and to
describe the relationship between singular (including regular)
Poisson reduction and by controllability distribution Poisson
reduction for CH system. This is our second goal of research in this
paper. The main contributions in this paper are given as follows.
(1) We give Poisson reducible conditions of CH system and symmetric
CH system by controllability distributions, and prove the Poisson
reducible property for CH systems to keep invariant under the
CH-equivalence; (2) We study the singular Poisson reduction and
SPR-CH-equivalence for CH system with symmetry, and prove the
singular Poisson reduction theorem, which is a generalization of the
regular Poisson reduction theorem in \cite{wazh12} to the singular
context; (3) We prove a theorem to explain the relationship between
Poisson reduction for singular Poisson reducible CH systems by
$G$-invariant controllability distribution and that for associated
reduced CH system by reduced controllability distribution.\\

A brief of outline of this paper is as follows. In the second
section, we review some relevant definitions and basic facts about
Poisson manifolds, generalized distributions, and the Poisson
reduction for Poisson manifold by distribution, as well as the CH
system defined by using Poisson tensor on the cotangent bundle of a
configuration manifold and its CH-equivalence, which will be used in
subsequent sections. The Poisson reduction of CH system by
controllability distribution and that of symmetric CH system by
$G$-invariant controllability distribution are considered
respectively, in the third section and the fourth section, and we
prove that the property of Poisson reduction for CH system by
controllability distribution leaves invariant under CH-equivalence
conditions. In the fifth section, the singular Poisson reducible CH
system and its SPR-CH-equivalence are considered, and singular
Poisson reduction theorem is proved, which shows the relationship
between the SPR-CH-equivalence for singular Poisson reducible CH
systems with symmetry and the CH-equivalence for associated singular
Poisson reduced CH systems. Moreover, in the sixth section, a
theorem is given to show the relationship between singular
(including regular) Poisson reduction and by controllability
distribution Poisson reduction for CH system. At last, some examples
are given to state theoretical results of Poisson reduction for CH
systems by controllability distributions. These research work
develop the reduction theory of controlled Hamiltonian systems with
symmetry and make us have much deeper understanding and recognition
for the structure of controlled Hamiltonian systems.

\section{Preliminaries}

In order to study the Poisson reductions of CH systems by
controllability distributions, we first give some relevant
definitions and basic facts about Poisson manifolds, generalized
distributions, and the Poisson reductions for Poisson manifolds by
distributions. We also recall briefly the CH systems defined by
using Poisson tensor on a Poisson fiber bundle and on the cotangent
bundle of a configuration manifold and their CH-equivalence, which
will be used in subsequent sections. we shall follow the notations
and conventions introduced in Abraham et al \cite{abma78,abmara88},
Marsden and Ratiu \cite{mara86, mara99}, Ortega and Ratiu
\cite{orra04}, Jotz and Ratiu \cite{jora09}, and Wang and Zhang
\cite{wazh12}.

\subsection{Poisson Manifolds and Generalized Distributions}

Let $P$ be a smooth manifold. $C^\infty(P)$ a set of smooth
functions on $P$, and $\mathfrak{X}(P)$ a set of smooth vector
fields on $P$. A Poisson bracket (or a Poisson structure) on the
manifold $P$ is a bilinear operation $\{\cdot,\cdot\}$ on
$C^\infty(P)$ such that $(C^\infty(P),\{\cdot,\cdot\})$ is a Lie
algebra; and $\{\cdot,\cdot\}$ is a derivation in each factor. A
manifold $P$ endowed with a Poisson bracket on $C^\infty(P)$ is
called a {\bf Poisson manifold}, and denoted by
$(P,\{\cdot,\cdot\})$. The pair $(C^\infty(P),\{\cdot,\cdot\})$ is
also called a {\bf Poisson algebra}. The derivation property of the
Poisson bracket implies that for any two functions $ f,g \in
C^\infty(P)$, the value of the bracket $\{f,g\}(m)$ at an arbitrary
point $m\in P$ depends on $f$ only through $\mathrm{d}f(m)$, which
allows us to define a covariant antisymmetric two-tensor $B\in
\Lambda^2(T^{\ast}P)$ by $B(m)(\mathrm{d}f(m),\mathrm{d}g(m))=\{f ,g
\}(m), \; \forall f,g \in C^\infty(P)$. This tensor $B$ is called
the {\bf Poisson tensor} of $P$. The vector bundle map $B^{\sharp}:
T^{\ast}P \rightarrow TP$ naturally associated to $B$ is defined by
$B(\alpha_m,\beta_m)= <\alpha_m, B^\sharp(m)\beta_m>, \; \forall
\alpha_m, \beta_m \in T_m^{\ast}P$. The range $D:=B^{\sharp}
(T^{\ast}P)\subset TP$ is called the {\bf characteristic
distribution}.\\

A {\bf generalized distribution} $D$ on $P$ is a subset of the
tangent bundle $TP$ such that for any point $m\in P$, the fiber
$D(m):=D\cap T_mP $ is a vector subspace of $T_mP $. The dimension
of $D(m)$ is called the {\bf rank} of $D$ at $m$. A point $m\in P$
is a regular point of the distribution $D$, if there exists a
neighborhood $U$ of $m$ such that the rank of $D$ is constant on
$U$. Otherwise, $m$ is a singular point of the distribution. A
distribution $D$ is called {\bf regular} if every point $m\in P$ is
a regular point of the distribution $D$. A {\bf differentiable
section} of $D$ is a differentiable vector field $X$ defined on an
open subset $U$ of $P$, such that for any point $z\in U$, $X(z)\in
D(z)$. An immersed connected submanifold $N$ of $P$ is said to be an
{\bf integral manifold} of the distribution $D$ if
$T_zi(T_zN)\subset D(z)$, $\forall z\in N$, where $i:N\to P$ is the
inclusion. $N$ is said to be of {\bf maximal dimension} at $z\in N$
if $T_zi(T_zN)=D(z)$. The generalized distribution $D$ is {\bf
differentiable} if for every point $m\in P$ and for every vector
$v\in D(m)$, there exist a differentiable section $X$ of $D$,
defined on an open neighborhood $U$ of $m$, such that $X(m)=v$. The
generalized distribution $D$ is {\bf completely integrable} if for
every point $m\in P$, there exists an integral manifold of $D$
everywhere of maximal dimension which contains $m$. The generalized
distribution $D$ is {\bf involutive} if it is invariant under the
(local) flows associated to differentiable sections of $D$. From
Stefen \cite{st74} and Sussmann \cite{su73} we know that $D$ is
completely integrable if and only if it is involutive.\\

Let $D$ be an integrable generalized distribution on $P$, then for
every point $m\in P$, there exists a unique connected integral
manifold $\mathcal{L}_m$ of $D$, which contains $m$ and has maximal
dimensions. $\mathcal{L}_m$ is called the {\bf maximal integral
manifold} or the accessible set of $D$ going through $m$, which is a
symplectic leaf of $P$. We know that the local structure of Poisson
manifolds is more complex than what one obtains in the symplectic
case. However, the {\bf symplectic stratification theorem} shows
that if $P$ is a finite dimensional Poisson manifold, then $P$ is
the disjoint union of its symplectic leaves. Each symplectic leaf in
$P$ is an injectively immersed Poisson submanifold and the induced
Poisson structure on the leaf is symplectic.\\

In the following we will be interested in the specific case in which
the generalized distribution is given by an everywhere defined
family of vector fields, that is, there is a family of smooth vector
fields $\mathcal{F}$ whose elements are vector fields $X$ defined on
a open subset Dom$(X)\subset P$ such that, for any $m\in P$ the
generalized distribution $D_{\mathcal{F}}$ is given by
$$D_{\mathcal{F}}(m) = \mbox{span}\{X(m)\in T_mP \mid X\in \mathcal{F}
\mbox{ and } m\in \mbox{Dom}(X)\}.$$ Note that in such a case the
distribution $D_{\mathcal{F}}$ is differentiable by construction. We
will say that $D_{\mathcal{F}}$ is the generalized distribution {\bf spanned} by $\mathcal{F}$.\\

Let $G$ be a Lie group acting properly and canonically on the
Poisson manifold $(P,\{\cdot,\cdot\})$ by the map $\Phi:G\times
P\rightarrow P $, and for any $g\in G$, the map
$\Phi_g:=\Phi(g,\cdot):P\rightarrow P$ is a diffeomorphism of $P$. A
submanifold $S (\subset P)$ is called $G$-{\bf invariant}, if
$\Phi_g(S)=S, \; \forall g \in G. $ A vector field $X \in
\mathfrak{X}(P)$ is called $G$-{\bf equivariant}, if $X \cdot \Phi_g
= T\Phi_g \cdot X, \; \forall g \in G. $ We denote the set of
$G$-equivariant vector fields on $P$ by $\mathfrak{X}(P)^G. $ Let
$\mathcal{F}$ be an everywhere defined family of local vector fields
on $P$ and $D_{\mathcal{F}}$ the generalized distribution spanned by
it. We say that the distribution $D_{\mathcal{F}}$ is $G$-{\bf
invariant} when for any $z\in P$ and any $g\in G, $ we have that
$T_z\Phi_g D_{\mathcal{F}}(z) = D_{\mathcal{F}}(\Phi_g(z)). $ A
distribution spanned by a family of $G$-equivariant vector fields is
always $G$-invariant, but the reverse implication is not necessarily
true, see Ortega and Ratiu \cite{orra04}. The following propositions
are important for the study in this paper.

\begin{prop}Let $(P,\{\cdot,\cdot\})$ be a Poisson manifold with Poisson tensor $B\in \Lambda^2(T^*P)$.
Then for any $m\in P$ and any vector subspace $V\subset T_mP $, we have that\\
{\bf (i)} $B^\sharp(m)(V^\circ)= (V\cap
T_m\mathcal{L})^{\omega_{\mathcal{L}(m)}};$\\
{\bf (ii)} $B^\sharp(m)(B^\sharp(m)(V^\circ))^\circ)= V\cap
T_m\mathcal{L},$\\
where $V^\circ:= \{\alpha_m \in T^*_mP |\; <\alpha_m, v>=0, \;
\forall v\in V \} \subset T^*P $ is the {\bf annihilator} of $V$ in
$T_m^*P$, and $\mathcal{L}_m $ is the {\bf symplectic leaf} of $P$
that contains the point $m$, and $\omega_{\mathcal{L}}$ is the
symplectic form on $\mathcal{L}$.
\end{prop}

\begin{prop}Let $(P,\{\cdot,\cdot\})$ be a Poisson manifold with Poisson tensor $B\in \Lambda^2(T^*P)$.
Assume that Lie group $G$ acts canonically on $P$, and for any $m\in
P$, $G_m$ is the {\bf isotropy subgroup} of $G$ at point $m$. Then we have that\\
{\bf (i)} $B^\sharp(m): T_m^*P \rightarrow T_mP $ is $G_m$-equivariant; \\
{\bf (ii)} If the Poisson bracket $\{\cdot,\cdot\}$ on $P$ is
induced by a symplectic form $\omega$, then for any vector subspace
$V\subset T_mP $, $B^\sharp(m)(V^{G_m})= (B^\sharp(m)(V))^{G_m}$,
where the $G_m$-superscript denotes the set of $G_m$-fixed points in
the corresponding space.
\end{prop}

In order to give some examples of application for the Poisson
reduction of CH systems, we need some specific submanifolds, for
examples, coisotropic submanifold, see Libermann and Marle \cite
{lima87}; and cosymplectic submanifold, see Weinstein \cite {we83}.

\begin{defi}
 Let $(P,\{\cdot,\cdot\})$ be a Poisson manifold with Poisson tensor $B\in
 \Lambda^2(T^*P)$, and $S\subset P$ an immersed smooth submanifold
 of $P$. The conormal bundle of the submanifold $S$ is
 $$(TS)^\circ:= \{ \alpha_s \in T^*_sP \; |\; <\alpha_s, v_s>=0, \; \forall s\in S,\; v_s \in T_sS \},$$
and it is a vector subbundle of $T^*P |_S. $\\
{\bf (i)} The submanifold $S$ is called {\bf coisotropic}, if
$B^\sharp((TS)^\circ) \subset TS; $\\
{\bf (ii)} The submanifold $S$ is called {\bf cosymplectic}, if $S$
is an embedded submanifold and satisfies that
$B^\sharp((TS)^\circ)\cap TS= \{0\}, $ and $T_sS +T_s\mathcal{L}_s =
T_sP, $ for any $s\in S$ and $\mathcal{L}_s $ the symplectic leaf of
$(P, B)$ containing $s\in S.$
\end{defi}

In particular, if $(P, B)$ is a symplectic manifold and its Poisson
bracket on $P$ is induced by a symplectic form $\omega \in
\Omega^2(P), $ then in this case the condition is given by
$$B^\sharp((TS)^\circ)= \{u\in TS \; |\; \omega(u, v)=0, \; \forall v \in (TS)^\circ \} = (TS)^{\omega} \subset TS, $$
that is, $S$ is a usually coisotropic submanifold of the symplectic
manifold $(P, \omega)$. Any cosymplectic submanifold of the
symplectic manifold $(P, \omega)$ is its symplectic submanifold.\\

\subsection{Reduction of Poisson Manifold by Distribution }

Let $(P,\{\cdot,\cdot\})$ be a Poisson manifold, and $S\subseteq P$
a decomposed subset of $P$. Let $\{S_i\}_{i\in I}$ be the pieces of
the corresponding decomposition. The topology of $S$ is not
necessarily the relative topology as a subset of $P$. We say that
$D\subset TP|_S$ is a {\bf smooth generalized distribution} (that
is, not necessarily of constant rank) on $S$ adapted to the
decomposition $\{S_i\}_{i\in I}$, if $D_{S_i}:=D\cap TS_i $ is a
smooth distribution on $S_i$ for all $i\in I$. The distribution $D$
is said to be {\bf integrable} if $D_{S_i}$ is integrable for each
$i\in I$. Thus, in this case, we can partition each $S_i$ into the
corresponding maximal integral manifolds, and the resulting
equivalence relations on each $S_i$, whose equivalence classes are
precisely these maximal integral manifolds, induce an equivalence
relation $D_S$ on the whole set $S$ by taking the union of the
different equivalence classes corresponding to all the $D_{S_i}$,
and we can define the quotient space $S/D_S:=\cup_{i\in
I}S_i/D_{S_i}$, and denote by $\pi_{D_S}:S\to S/D_S$ the natural
projection. We say that a subset $U \subset S$ is {\bf
$D_S$-invariant} if it is invariant under the flow
of any section of $D_{S_i}$, for all $i\in I$.\\

We define the {\bf presheaf} of smooth functions $C_{S/D_S}^\infty$
on $S/D_S$ as the map that associates to any open subset $V$ of
$S/D_S$ the set of functions $C_{S/D_S}^\infty(V)$ characterized by
the following property: $f\in C_{S/D_S}^\infty(V)$ if and only if
for any $z\in V$, there exists $m\in \pi_{D_S}^{-1}(z)$, $U_m$ an
open neighborhood of $m\in P$, such that $U_m\cap S$ is
$D_S$-invariant, and $F\in C_P^\infty(U_m)$ satisfying
\begin{equation}
f\circ \pi_{D_S}|_{\pi_{D_S}^{-1}(V)\cap
U_m}=F|_{\pi_{D_S}^{-1}(V)\cap U_m}.\label{localex}
\end{equation}
In this case, we say that $F$ is a {\bf local extension} of $f\circ
\pi_{D_S}$ at the point $m$. Moreover, the presheaf
$C_{S/D_S}^\infty$ is said to have the {\bf $(D,D_S)$-local
extension property}, when the topology of $S$ is stronger than the
relative topology and the local extensions of $f\circ \pi_{D_S}$
defined in (\ref{localex}) can always be chosen to satisfy
\begin{equation*}
    \mathrm{d}F(u)|_{D(u)}=0,\; \mbox{for any } u\in
    \pi_{D_S}^{-1}(V)\cap U_m.
\end{equation*}
In this case, we say that $F$ is a {\bf local $D$-invariant
extension} of $f\circ \pi_{D_S}$ at the point $m$.\\

At first, we recall the following definition, see Ortega and Ratiu
\cite{orra04}.
\begin{defi}
  Let $(P,\{\cdot,\cdot\})$ be a Poisson manifold, $S\subset P$ a
  decomposed space, and $D\subset TP|_S$ a smooth distribution
  adapted to the decomposition $\{S_i\}_{i\in I}$ of $S$. The
  distribution $D$ is called {\bf Poisson or canonical}, if the condition
  $\mathrm{d}f|_D=\mathrm{d}g|_D=0$, for any $f,g\in C_P^\infty(U)$
  and any open subset $U\subset P$, implies that
  $\mathrm{d}\{f,g\}|_D=0$.
\end{defi}

From this definition, it is easy to know immediately the following
result. Assume that distribution $D$ is Poisson and the presheaf
$C_{S/D_S}^\infty$ has the $(D,D_S)$-local extension property. If
$F$ and $H \in C_P^\infty(U_m)$ are local $D$-invariant extensions
of $f\circ \pi_{D_S}$ and $h \circ \pi_{D_S}$ at $m\in
\pi_{D_S}^{-1}(z)$, then their Poisson bracket $\{F,H\}$ is also
local $D$-invariant extension.\\

Now we can introduce the definition of reduction by distribution for
Poisson manifold as follows.
\begin{defi}
  Let $(P,\{\cdot,\cdot\})$ be a Poisson manifold with Poisson tensor $B\in \Lambda^2(T^*P)$,
  $S$ a decomposed subset of $P$, and $D\subset TP|_S$ a Poisson integrable
  distribution that is adapted to the decomposition of $S$, and $D_{S}:=D\cap TS $.
  Assume that the presheaf
  $C_{S/D_S}^\infty$ has the $(D,D_S)$-local extension
  property. We say that $(P,B,D,S)$ is {\bf Poisson
  reducible}, when $(S/D_S,C_{S/D_S}^\infty,B^{S/D_S})$
  is a well-defined presheaf of Poisson algebras, where for any open
  set $V\subset S/D_S$, the Poisson tensor $B^{S/D_S}$ is defined by the Poisson bracket
  $$\{\cdot,\cdot\}_V^{S/D_S}: C_{S/D_S}^\infty(V)\times
  C_{S/D_S}^\infty(V)\to C_{S/D_S}^\infty(V),$$
 given by $$\{f,h\}_V^{S/D_S}(\pi_{D_S}(m)):=\{F,H\}_B(m),$$ for any
  $m\in \pi_{D_S}^{-1}(z), \, z\in V$ and any local $D$-invariant extension $F$ and
  $H \in C_{P}^\infty(U_m)$ of $f\circ \pi_{D_S}$ and $h \circ \pi_{D_S}$ at $m$,
  respectively, where $\pi_{D_S}: S \rightarrow S/D_S$ is the projection, and $U_m$ is
an open neighborhood of $m$, such that $U_m\cap S$ is
$D_S$-invariant. That is, $B^{S/D_S}(\pi_{D_S}(m)):= B(m).$
\end{defi}

The following theorem give the Poisson reducible conditions by
distribution for Poisson manifold, see Ortega and Ratiu
\cite{orra04} and Jotz and Ratiu \cite {jora09}.
\begin{theo}\label{poissonreducible}
  Assume that $(P,\{\cdot,\cdot\})$ is a Poisson manifold with
  Poisson tensor $B\in \Lambda^2(T^*P)$, $S$ a
  decomposed set of $P$, and $D\subset TP|_S$ a Poisson integrable
  generalized distribution adapted to the decomposition of $S$,
  $D_{S}:=D\cap TS $,
  such that the presheaf $C_{S/D_S}^\infty$ has the $(D,D_S)$-local extension
  property. Then $(P,B,D,S)$ is Poisson reducible if
  and only if for any $m\in S$
  \begin{equation}
    B^\sharp(\Delta_m)\subset (\Delta_m)_S^\circ\label{prdcon}
  \end{equation}
  where \begin{equation*}\Delta_m:=\left\{\mathrm{d}F(m)\left|\begin{aligned}F\in
  C_P^\infty(U_m),\; \mathrm{d}F(s)|_{D(s)}=0\; \mbox{for all } s\in
  U _m\cap S,\\ \mbox{ and for any open neighborhood }U_m \mbox{ of } m \mbox{ in }
  P\end{aligned}\right.\right\},\end{equation*}
  and \begin{equation*}
    (\Delta_m)_S:=\left\{\mathrm{d}F(m)\in \Delta_m\left|\begin{aligned}F|_{U_m\cap V_m}\
    \mbox{is constant for any open neighborhood } U_m\\
    \mbox{ of } m \mbox{ in } P \mbox{ and an open neighborhood } V_m \mbox{ of } m \mbox{ in }
    S\end{aligned}\right.\right\}.
  \end{equation*}
In particular, if $S$ is endowed with the relative topology, then
the definition of $(\Delta_m)_S$ simplifies to \begin{equation*}
    (\Delta_m)_S=\{\mathrm{d}F(m)\in \Delta_m|\;F|_{U_m\cap S}
    \mbox{ is constant for an open neighborhood } U_m
    \mbox{ of } m \mbox{ in } P\}.
  \end{equation*}
\end{theo}

When $S$ is an embedded submanifold of $P$ and $D \subset TP|_S$ is
a subbundle of the tangent bundle of $P$ restricted to $S$. In this
case Ortega and Ratiu in \cite{orra04} gave the following
proposition.
\begin{prop}
Assume that $(P,\{\cdot,\cdot\})$ is a Poisson manifold and $S$ is
an embedded submanifold of $P$ and $D \subset TP|_S$ is a subbundle
of the tangent bundle of $P$ restricted to $S$, such that $D_S=D\cap
TS$ is a smooth regular integrable distribution on $S$. Then the
presheaf $C_{S/D_S}^\infty$ has the $(D,D_S)$-local extension
property.
\end{prop}

Moreover, by using Theorem 2.6 and the above proposition, one can
also obtain another result on the Poisson reducible conditions by
distribution for Poisson manifold, which are given in Falceto and
Zambon \cite{faza08} and Marsden and Ratiu \cite{mara86}. In
particular, for some specific submanifolds, including the
coisotropic and cosymplectic submanifolds, by using Proposition 2.1
one can get the following reduction theorem, see Ortega and Ratiu
\cite {orra04}.

\begin{theo}
Let $(P,\{\cdot,\cdot\})$ be a Poisson manifold with Poisson tensor
$B\in \Lambda^2(T^*P)$, and $S\subset P$ an embedded smooth
submanifold of $P$ and the distribution $D:= B^\sharp((TS)^\circ)
\subset TP|_S.$ If any one of the following three conditions
holds,\\
{\bf (i)} $D_S:= D\cap TS$ is a smooth and integrable generalized
distribution on $S$, such that the presheaf $C_{S/D_S}^\infty$ has
the $(D,D_S)$-local extension property;\\
{\bf (ii)} $S$ is coisotropic, and the presheaf $C_{S/D_S}^\infty$
has the
$(D,D_S)$-local extension property;\\
{\bf (iii)} $S$ is cosymplectic;\\
Then $(P, B, D, S)$ is Poisson reducible.
\end{theo}

\subsection{Controlled Hamiltonian System and CH-equivalence}

In this subsection, we shall introduce the CH system defined by
using Poisson structure on the cotangent bundle of a configuration
manifold, which is a special case of the Definition 3.1 of CH system
defined on a Poisson fiber bundle in \cite{wazh12}. We also discuss
the controlled Hamiltonian equivalence (CH-equivalence) of such
systems. For convenience, we assume that all
controls appearing in this paper are the admissible controls.\\

Let $(E,M,N,\pi,G)$ be a Poisson fiber bundle and $B$ be a Poisson
tensor on $E$, then we have a induced bundle map $B^{\sharp}:
T^{\ast}E \rightarrow TE $ such that for any $\lambda_z, \; \nu_z
\in T^{\ast}_zE, \; z \in E, \; B(\lambda_z,\nu_z)= <\lambda_z,
B^\sharp(z)\nu_z>$. If $H: E \rightarrow \mathbb{R}$ is a
Hamiltonian, then the Hamiltonian vector field $X_H \in TE $ can be
expressed by $X_H = B^\sharp \mathbf{d}H$, and $(E, B, H )$ is a
Hamiltonian system. Moreover, if considering the external force and
control, we can define a CH system on the Poisson fiber bundle $E$
as follows.

\begin{defi}
(CH System) A {\bf CH system} on $E$ is a 5-tuple $(E, B, H, F, W)$,
where $(E, B, H )$ is a Hamiltonian system, and the function $H: E
\rightarrow \mathbb{R}$ is called the Hamiltonian, a
fiber-preserving map $F: E\rightarrow E$ is called the (external)
force map, and a fiber submanifold $W$ of $E$ is called the control
subset.
\end{defi}
Sometimes, $W$ also denotes the set of fiber-preserving maps from
$E$ to $W$. When a feedback control law $u: E \rightarrow W$ is
chosen, the 5-tuple $(E, B, H, F, u)$ denotes a closed-loop dynamic
system.\\

Let $Q$ be a smooth manifold, and $T^\ast Q$ its cotangent bundle.
If $T^\ast Q$ has a Poisson structure $\{\cdot,\cdot\}$, we can
define a Poisson tensor $B$ on $T^{\ast}Q$ such that $(T^\ast Q, B)$
is a Poisson vector bundle. If we take that $E= T^* Q$, from above
definition we can obtain a CH system on the cotangent bundle $T^\ast
Q$, that is, 5-tuple $(T^\ast Q, B, H, F, W)$. Where the
fiber-preserving map $F: T^*Q\rightarrow T^*Q$ is the (external)
force map, that is the reason that the fiber-preserving map $F:
E\rightarrow E$ is called an (external) force map in above
definition. In particular, the cotangent bundle $T^\ast Q$ has a
canonical symplectic form $\omega $, so $(T^\ast Q, \omega )$ is a
symplectic vector bundle. From above definition we also obtain a RCH
system on the cotangent bundle $T^\ast Q$, that is, 5-tuple $(T^\ast
Q, \omega, H, F, W)$, see Marsden et al \cite{mawazh10}.\\

In order to describe the dynamics of the CH system $(E,B,H,F,W)$
with a control law $u$, we need to give a good expression of the
dynamical vector field of CH system. At first, we introduce a
notations of vertical lift maps of a vector along a fiber. For a
smooth manifold $E$, its tangent bundle $TE$ is a vector bundle, and
for the fiber bundle $\pi: E \rightarrow M$, we consider the tangent
mapping $T\pi: TE \rightarrow TM$ and its kernel $ker
(T\pi)=\{\rho\in TE| T\pi(\rho)=0\}$, which is a vector subbundle of
$TE$. Denote by $VE:= ker(T\pi)$, which is called a {\bf vertical
bundle} of $E$. Assume that there is a metric on $E$, and we take a
Levi-Civita connection $\mathcal{A}$ on $TE$, and denote by $HE:=
ker(\mathcal{A})$, which is called a {\bf horizontal bundle} of $E$,
such that $TE= HE \oplus VE. $ For any $x\in M, \; a_x, b_x \in E_x,
$ any tangent vector $\rho(b_x)\in T_{b_x}E$ can be split into
horizontal and vertical parts, that is, $\rho(b_x)=
\rho^h(b_x)\oplus \rho^v(b_x)$, where $\rho^h(b_x)\in H_{b_x}E$ and
$\rho^v(b_x)\in V_{b_x}E$. Let $\gamma$ be a geodesic in $E_x$
connecting $a_x$ and $b_x$, and denotes by $\rho^v_\gamma(a_x)$ a
tangent vector at $a_x$, which is a parallel displacement of the
vertical vector $\rho^v(b_x)$ along the geodesic $\gamma$ from $b_x$
to $a_x$. Since the angle between two vectors is invariant under a
parallel displacement along a geodesic, then
$T\pi(\rho^v_\gamma(a_x))=0, $ and hence $\rho^v_\gamma(a_x) \in
V_{a_x}E. $ Now, for $a_x, b_x \in E_x $ and tangent vector
$\rho(b_x)\in T_{b_x}E$, we can define the {\bf vertical lift map}
of a vector along a fiber given by
$$\mbox{vlift}: TE_x \times E_x \rightarrow TE_x; \;\; \mbox{vlift}(\rho(b_x),a_x) = \rho^v_\gamma(a_x). $$
It is easy to check from the basic fact in differential geometry
that this map does not depend on the choice of $\gamma$. If $F: E
\rightarrow E$ is a fiber-preserving map, for any $x\in M$, we have
that $F_x: E_x \rightarrow E_x$ and $TF_x: TE_x \rightarrow TE_x$,
then for any $a_x \in E_x$ and $\rho\in TE_x$, the vertical lift of
$\rho$ under the action of $F$ along a fiber is defined by
$$(\mbox{vlift}(F_x)\rho)(a_x)=\mbox{vlift}((TF_x\rho)(F_x(a_x)), a_x)= (TF_x\rho)^v_\gamma(a_x), $$
where $\gamma$ is a geodesic in $E_x$ connecting $F_x(a_x)$ and
$a_x$.\\

In particular, when $\pi: E \rightarrow M$ is a vector bundle, for
any $x\in M$, the fiber $E_x=\pi^{-1}(x)$ is a vector space. In this
case, we can choose the geodesic $\gamma$ to be a straight line, and
the vertical vector is invariant under a parallel displacement along
a straight line, that is, $\rho^v_\gamma(a_x)= \rho^v(b_x).$
Moreover, when $E= T^*Q, \; M=Q $, by using the local trivialization
of $TT^*Q$, we have that $TT^*Q\cong TQ \times T^*Q$. Because of
$\pi: T^*Q \rightarrow Q$, and $T\pi: TT^*Q \rightarrow TQ$, then in
this case, for any $\alpha_x, \; \beta_x \in T^*_x Q, \; x\in Q, $
we know that $(0, \beta_x) \in V_{\beta_x}T^*_x Q, $ and hence we
can get that
$$ \mbox{vlift}((0, \beta_x)(\beta_x), \alpha_x) = (0, \beta_x)(\alpha_x)
=\left.\frac{\mathrm{d}}{\mathrm{d}s}\right|_{s=0}(\alpha_x+s\beta_x),
$$ which is consistent with the definition of vertical lift map
along
fiber in Marsden and Ratiu \cite{mara99}.\\

For a given CH System $(T^\ast Q, B, H, F, W)$, the dynamical vector
field of the associated Hamiltonian system $(T^\ast Q, B, H) $ is
that $X_H= B^\sharp\mathrm{d}H$, where, $B^\sharp: T^\ast T^\ast Q
\rightarrow TT^\ast Q; \mathrm{d}H\mapsto B^\sharp\mathrm{d}H, $
such that for any $\lambda \in T^\ast T^\ast Q$,
$B(\lambda,\mathrm{d}H)= <\lambda, B^\sharp\mathrm{d}H>$. If
considering the external force $F: T^*Q \rightarrow T^*Q, $ by using
the above notations of vertical lift maps of a vector along a fiber,
the change of $X_H$ under the action of $F$ is that
$$\mbox{vlift}(F)X_H(\alpha_x)= \mbox{vlift}((TFX_H)(F(\alpha_x)), \alpha_x)= (TFX_H)^v_\gamma(\alpha_x),$$
where $\alpha_x \in T^*_x Q, \; x\in Q $ and $\gamma$ is a straight
line in $T^*_x Q$ connecting $F_x(\alpha_x)$ and $\alpha_x$. In the
same way, when a feedback control law $u: T^\ast Q \rightarrow W$ is
chosen, the change of $X_H$ under the action of $u$ is that
$$\mbox{vlift}(u)X_H(\alpha_x)= \mbox{vlift}((TuX_H)(u(\alpha_x)), \alpha_x)= (TuX_H)^v_\gamma(\alpha_x).$$
In consequence, we can give an expression of the dynamical vector
field of CH system as follows.

\begin{prop}
The dynamical vector field of a CH system $(T^\ast Q,B,H,F,W)$ with
a control law $u$ is the synthetic of Hamiltonian vector field $X_H$
and its changes under the actions of the external force $F$ and
control $u$, that is,
$$X_{(T^\ast Q,B,H,F,u)}(\alpha_x)= X_H(\alpha_x)+ \textnormal{vlift}(F)X_H(\alpha_x)
+ \textnormal{vlift}(u)X_H(\alpha_x),$$
for any $\alpha_x \in T^*_x Q, \; x\in Q $. For convenience, it is
simply written as
\begin{equation}X_{(T^\ast Q,B,H,F,u)}
=B^\sharp\mathrm{d}H +\textnormal{vlift}(F) +\textnormal{vlift}(u).
\end{equation}
\end{prop}
We also denote that $\mbox{vlift}(W)=\bigcup \{\mbox{vlift}(u)X_H |
\; u\in W\}$. For the CH system $(E,B,H,F,W)$ with a control law
$u$, we have also a similar expression of its dynamical vector
field. It is worthy of note that in order to deduce and calculate
easily, we always use the simple expression of dynamical vector
field $X_{(T^\ast
Q,B,H,F,u)}$.\\

Next, we note that when a CH system is given, the force map $F$ is
determined, but the feedback control law $u: T^\ast Q\rightarrow W$
could be chosen. In order to describe the feedback control law to
modify the structure of CH system, we give the controlled
Hamiltonian matching conditions and CH-equivalence as follows.
\begin{defi}
(CH-equivalence) Suppose that we have two CH systems $(T^\ast
Q_i,B_i,H_i,F_i,W_i)$, $i=1,2$, we say them to be {\bf
CH-equivalent}, or simply, $(T^\ast
Q_1,B_1,H_1,F_1,W_1)\stackrel{CH}{\sim}(T^\ast
Q_2,B_2,H_2,\\F_2,W_2)$, if there exists a diffeomorphism $\varphi:
Q_1\rightarrow Q_2$, such that the following Hamiltonian matching
conditions hold:

\noindent {\bf HM-1:} The cotangent lift map of $\varphi$, that is,
$\varphi^\ast= T^\ast \varphi: T^\ast Q_2\rightarrow T^\ast Q_1$ is
a Poisson map, and $W_1=\varphi^\ast (W_2)$;

\noindent {\bf HM-2:}
$Im[B_1^\sharp\mathrm{d}H_1+\textnormal{vlift}(F_1)-T\varphi^\ast
(B_2^\sharp\mathrm{d}H_2)-\textnormal{vlift}(\varphi^\ast
F_2\varphi_\ast)]\subset \textnormal{vlift}(W_1)$, where the map
$\varphi_\ast=(\varphi^{-1})^\ast: T^\ast Q_1\rightarrow T^\ast
Q_2$, $T\varphi^\ast: TT^\ast Q_2 \rightarrow TT^\ast Q_1$ is the
tangent map of $\varphi^\ast$ and $Im$ means the pointwise image of
the map in brackets.
\end{defi}

It is worthy of note that our CH-system is defined by using the
Poisson tensor on the cotangent bundle of a configuration manifold,
we must keep with the Poisson structure when we define the
CH-equivalence, that is, the induced equivalent map $\varphi^*$ is
Poisson on the cotangent bundle. The following Theorem 2.12 explains
the significance of the above CH-equivalence relation, the proof is
given in \cite {wazh12}.
\begin{theo}
Suppose that two CH systems $(T^\ast Q_i,B_i,H_i,F_i,W_i)$, $i=1,2,$
are CH-equivalent, then there exist two control laws $u_i: T^\ast
Q_i \rightarrow W_i, \; i=1,2, $ such that the two associated
closed-loop systems produce the same equations of motion, that is,
$X_{(T^\ast Q_1,B_1,H_1,F_1,u_1)}\cdot \varphi^\ast = T\varphi^\ast
 X_{(T^\ast Q_2,B_2,H_2,F_2,u_2)}$. Moreover, the explicit relation between
the two control laws $u_i, i=1,2,$ is given by
\begin{equation}\textnormal{vlift}(u_1)- \textnormal{vlift}(\varphi^\ast u_2\varphi_\ast)
=-B_1^\sharp\mathrm{d}H_1-\textnormal{vlift}(F_1)+T\varphi^\ast
(B_2^\sharp\mathrm{d}H_2)+\textnormal{vlift}(\varphi^\ast F_2
\varphi_\ast).
\end{equation}
\end{theo}

\section{Poisson Reduction of CH System by Controllability Distribution}

In this section, we first give a definition of controllability
distribution of CH system, then give a theorem to show the Poisson
reducible conditions of CH system by controllability distribution.
Moreover, we prove that the property of Poisson reduction for CH
system by controllability distribution leaves invariant under the
CH-equivalence.

\begin{defi}
For a CH system $(T^*Q, B, H, F, W)$, the fiber submanifold $W$ of
$T^*Q$ is its control subset. If for each $u \in W, $ the dynamical
control system $\dot{x}= X_{(T^\ast Q,B,H,F,u)}$ is controllable,
that is, for any two states $x_0$ and $x_1$ of this system, there is
a finite piecewise smooth integral curve $x(t)$ of $\dot{x}=
X_{(T^\ast Q,B,H,F,u)}, \; t\in [0,1]$, such that $x(0)=x_0$ and
$x(1)=x_1$. Then $W$ is called a {\bf controllability submanifold}
of CH system.
\end{defi}
\begin{defi}
Assume that $W$ is a controllability submanifold of CH system
$(T^*Q, B, H, F, W)$, and a distribution $D \subset TT^*Q|_W $, and
it is a Poisson integrable generalized distribution, then $D$ is
called a {\bf controllability distribution} of CH system.
\end{defi}

For a CH system $(T^*Q, B, H, F, W)$, assume that $W$ is its
controllability submanifold and $D \subset TT^*Q|_W $ is a
controllability distribution such that $D_W=D\cap TW$ is a smooth
regular integrable distribution on $W$. If the presheaf
$C_{W/D_W}^\infty$ has the $(D,D_W)$-local extension property, then
we may consider that $(T^*Q,B,D,W)$ is Poisson reducible if
$(W/D_W,C_{W/D_W}^\infty,B^{W/D_W})$ is a well-defined presheaf of
Poisson algebras. Thus, by using the Poisson reduction by
controllability distribution $D$ for the phase space of CH system,
we can define the Poisson reducible CH system by controllability
distribution as follows.
\begin{defi}
A CH system $(T^*Q, B, H, F, W)$ is called to be {\bf Poisson
reducible by controllability distribution} $D$, if $(T^*Q,B,D,W)$ is
Poisson reducible in the sense of Definition 2.5.
\end{defi}
From above definition and Theorem 2.6, we can obtain the following
theorem.
\begin{theo}
Suppose that $(T^*Q, B, H, F, W)$ is a CH system, and $W$ is its
controllability submanifold and $D\subset TT^*Q|_W$ is its
controllability distribution, such that $D_W=D\cap TW$ is a smooth
regular integrable distribution on $W$, and the presheaf
$C_{W/D_W}^\infty$ has the $(D,D_W)$-local extension property. Then
the CH system $(T^*Q, B, H, F, W)$ is Poisson reducible by
controllability distribution $D$, if and only if for any $z \in W$
  \begin{equation}
    B^\sharp(\Delta_z)\subset (\Delta_z)_W^\circ\label{prdcon}
  \end{equation}
where
\begin{equation*}\Delta_z:=\left\{\mathrm{d}f(z)\left|\begin{aligned}f\in
  C_{T^*Q}^\infty(U_z),\; \mathrm{d}f(y)|_{D(y)}=0,\; \mbox{for all } y\in
  U _z\cap W,\\ \mbox{ and for any open neighborhood }U_z \mbox{ of } z \mbox{ in }
  T^*Q\end{aligned}\right.\right\},\end{equation*}
  and \begin{equation*}
    (\Delta_z)_W:=\left\{\mathrm{d}f(z)\in \Delta_z\left|\begin{aligned}f|_{U_z\cap W}\
    \mbox{is constant for any open}\\ \mbox{ neighborhood } U_z
    \mbox{ of } z \mbox{ in } T^*Q \end{aligned}\right.\right\}.
  \end{equation*}
\end{theo}
{\bf Proof.} In fact, we take that $P=T^*Q $ and $S= W$, the
conclusion is a direct consequence of Theorem 2.6. $\hskip 1cm
\blacksquare$

It is worthy of note that for convenience, here and in subsequent
sections $W$ is endowed with the relative topology. For the general
case that the topology of $W$ is stronger than the relative
topology, we can also obtain the similar result. Moreover, if
considering the CH-equivalence of CH system, we can get the
following theorem to state that the property of Poisson reduction
for CH system by controllability distribution leaves invariant under
CH-equivalence conditions.

\begin{theo} Suppose that two CH systems
$(T^\ast Q_i,B_i,H_i,F_i,W_i)$, $i=1,2,$ are CH-equivalent with
equivalent map $\varphi^*: T^\ast
Q_2\rightarrow T^\ast Q_1$. Then we have that\\

\noindent{\bf (i)} $W_1$ is controllability submanifold of CH system
$(T^*Q_1, B_1, H_1, F_1, W_1)$ if and only if $W_2$ is
controllability submanifold
of CH system $(T^*Q_2, B_2, H_2, F_2, W_2). $\\

\noindent{\bf (ii)} $D_1=T\varphi^*(D_2) \subset TT^*Q_1|_{W_1}$ is
controllability distribution of CH system $(T^*Q_1, B_1, H_1, F_1,
W_1), $ such that $D_{W_1}=D_1\cap TW_1$ is smooth regular
integrable distributions on $W_1$ and the presheaf
$C_{W_1/D_{W_1}}^\infty$ has the $(D_1,D_{W_1})$-local extension
property if and only if $D_2 \subset TT^*Q_2|_{W_2}$ is
controllability distribution of CH system $(T^*Q_2, B_2, H_2, F_2,
W_2), $ such that $D_{W_2}=D_2\cap TW_2$ is smooth regular
integrable distributions on $W_2$ and the presheaf
$C_{W_2/D_{W_2}}^\infty$ has the $(D_2,D_{W_2})$-local extension
property.\\

\noindent{\bf (iii)} CH system $(T^*Q_1, B_1, H_1, F_1, W_1)$ is
Poisson reducible by controllability distribution $D_1$ if and only
if CH system $(T^*Q_2, B_2, H_2, F_2, W_2)$ is Poisson reducible by
controllability distribution $D_2$.
\end{theo}

\noindent{\bf Proof.} {\bf (i)} Because CH systems $(T^\ast
Q_i,B_i,H_i,F_i,W_i)$, $i=1,2,$ are CH-equivalent, then there exists
a diffeomorphism $\varphi: Q_1\rightarrow Q_2$, such that the
cotangent lift map $\varphi^\ast= T^\ast \varphi: T^\ast
Q_2\rightarrow T^\ast Q_1$ and the tangent lift map $T\varphi^\ast:
TT^\ast Q_2 \rightarrow TT^\ast Q_1$ are vector bundle isomorphism,
and by Theorem 2.12 there exist always a pair of control laws $u_i:
T^\ast Q_i \rightarrow W_i, \; i=1,2, $ such that $X_{(T^\ast
Q_1,B_1,H_1,F_1,u_1)}\cdot \varphi^\ast = T\varphi^\ast
 X_{(T^\ast Q_2,B_2,H_2,F_2,u_2)}$. Note that $W_1=\varphi^\ast
 (W_2)$, hence from Definition 3.1 we know that {\bf (i)} holds.\\

\noindent{\bf (ii)} Notice that $D_1=T\varphi^*(D_2), $ and $D_2
\subset TT^*Q_2|_{W_2}$. Thus,
$$D_{W_1}=D_1\cap TW_1 = T\varphi^*(D_2) \cap T\varphi^\ast
(TW_2)= T\varphi^\ast (D_2 \cap TW_2)= T\varphi^\ast (D_{W_2}).$$
Since $\varphi^\ast$ and $T\varphi^\ast$ are vector bundle
isomorphism, and $\varphi^\ast$ is Poisson, from Definition 2.4,
Definition 3.2 and the definition of the presheaf local extension
property, we know that {\bf (ii)} holds.\\

\noindent{\bf (iii)} If CH system $(T^*Q_1, B_1, H_1, F_1, W_1)$ is
Poisson reducible by controllability distribution $D_1$, we shall
prove that the CH system $(T^*Q_2, B_2, H_2, F_2, W_2)$ is Poisson
reducible by controllability distribution $D_2$. By the above
conclusions {\bf (i), (ii)} and Theorem 3.4, it suffices to show
that for any $z_2 \in W_2, $ we have that
$B_2^\sharp(\Delta^2_{z_2})\subset (\Delta^2_{z_2})_{W_2}^\circ $.
In fact, assume that for any $\alpha_2=\mathrm{d}f_2(z_2)\in
\Delta^2_{z_2}, $ that is, $f_2 \in C_{T^*Q_2}^\infty(U^2_{z_2}), $
where $U^2_{z_2}$ is an open neighborhood of $z_2$ in $T^*Q_2, $
such that $U^2_{z_2}\cap W_2$ is $D_{W_2}$-invariant and $
\mathrm{d}f_2(s_2)|_{D_2(s_2)}=0,\; \mbox{for all } s_2\in
  U^2_{z_2}\cap W_2. $ Since two CH systems $(T^\ast Q_i,B_i,H_i,F_i,W_i)$, $i=1,2,$
are CH-equivalent, and the equivalent map $\varphi^*: T^*Q_2
\rightarrow T^*Q_1$ is a Poisson diffeomorphism, then for $z_1 \in
T^*Q_1, $ such that $z_2= \varphi^*(z_1), $ there is an open
neighborhood $U^1_{z_1}$ of $z_1$ in $T^*Q_1, $ such that
$U^2_{z_2}= \varphi^*(U^1_{z_1})$ and $U^1_{z_1}\cap W_1$ is
$D_{W_1}$-invariant, and there exists $f_1 \in
C_{T^*Q_1}^\infty(U^1_{z_1}), $ such that $f_2= f_1\cdot \varphi^*$
and $ \mathrm{d}f_1(s_1)|_{D_1(s_1)}=0,\; \mbox{for all } s_1\in
  U^1_{z_1}\cap W_1, $ and
$\alpha_1=\mathrm{d}f_1(z_1)\in \Delta^1_{z_1}. $ In the same way,
for any $\beta_2=\mathrm{d}g_2(z_2)\in (\Delta^2_{z_2})_{W_2}, $
that is, $g_2 \in C_{T^*Q_2}^\infty(U^2_{z_2}) $ and
$g_2|_{U^2_{z_2}\cap W_2}$ is constant for any open neighborhood
$U^2_{z_2}$ of $z_2$ in $T^*Q_2$, then there is a $g_1 \in
C_{T^*Q_1}^\infty(U^1_{z_1}), $ such that $g_2= g_1\cdot \varphi^*$
and $g_1|_{U^1_{z_1}\cap W_1}$ is constant for the corresponding
open neighborhood $U^1_{z_1}$ of $z_1$ in $T^*Q_1$, and
$\beta_1=\mathrm{d}g_1(z_1)\in (\Delta^1_{z_1})_{W_1}. $ Because CH
system $(T^*Q_1, B_1, H_1, F_1, W_1)$ is Poisson reducible by
controllability distribution $D_1$, from the above conclusions {\bf
(i), (ii)} and Theorem 3.4 we have that
$B_1^\sharp(\Delta^1_{z_1})\subset (\Delta^1_{z_1})_{W_1}^\circ. $
It follows that
$$\{f_1,g_1\}_{B_1}(z_1)= <\mathrm{d}g_1(z_1),
B_1^\sharp(\mathrm{d}f_1(z_1))> = <\beta_1, B_1^\sharp(\alpha_1)>
=0.$$ Notice that the map $\varphi^*: T^*Q_2 \rightarrow T^*Q_1$ is
Poisson, we have that
\begin{align*}<\beta_2, B_2^\sharp(\alpha_2)> &=
<\mathrm{d}g_2(z_2), B_2^\sharp(\mathrm{d}f_2(z_2))>=
\{f_2,g_2\}_{B_2}(z_2)\\ &= \{f_1\cdot \varphi^*,g_1\cdot
\varphi^*\}_{B_2}(z_2)= (\varphi^*)^* \{f_1,g_1\}_{B_1}(z_1)=0.
\end{align*} Thus, we prove that $B_2^\sharp(\alpha^2)\in (\Delta^2_{z_2})_{W_2}^\circ$
and $B_2^\sharp(\Delta^2_{z_2})\subset (\Delta^2_{z_2})_{W_2}^\circ
$, and hence CH system\\ $(T^*Q_2, B_2, H_2, F_2, W_2)$ is Poisson
reducible by controllability distribution $D_2$.\\

Conversely, if CH system $(T^*Q_2, B_2, H_2, F_2, W_2)$ is Poisson
reducible by controllability distribution $D_2$, by using the same
way we can verify that for any $z_1 \in W_1, $ we have that
$B_1^\sharp(\Delta^1_{z_1})\subset (\Delta^1_{z_1})_{W_1}^\circ $.
Thus, the CH system $(T^*Q_1, B_1, H_1, F_1, W_1)$ is Poisson
reducible by controllability distribution $D_1$ by the above
conclusions {\bf (i), (ii)} and Theorem 3.4. $\hskip 1cm
\blacksquare$

\section{Poisson Reduction of Symmetric
CH System by $G$-invariant Controllability Distribution}

In this section, we shall consider the symmetric CH system and give
the Poisson reducible conditions of this system by $G$-invariant
controllability distribution.\\

Let $Q$ be a smooth manifold and $T^\ast Q$ its cotangent bundle
with a Poisson tensor $B$. Let $\Phi: G\times Q\rightarrow Q$ be a
smooth left action of the Lie group $G$ on $Q$, and the cotangent
lifted left action $\Phi^{T^\ast}:G\times T^\ast Q\rightarrow T^\ast
Q$ is Poisson. If Hamiltonian $H: T^{\ast}Q \rightarrow \mathbb{R}$
is $G$-invariant, then the 4-tuple $(T^{\ast}Q, G, B, H )$ is a
symmetric Hamiltonian system. Moreover, if the external force map
$F: T^\ast Q\rightarrow T^\ast Q$ and control subset $W$ of
\;$T^\ast Q$ are $G$-invariant, then we can give the definition of
symmetric CH system on $T^\ast Q$ as follows.
\begin{defi}
(Symmetric CH System) A {\bf symmetric CH system} on $T^{\ast}Q$ is
a 6-tuple $(T^\ast Q, G, B, H, F, W)$, where the action of $G$ on
$T^\ast Q$ is Poisson, and the Hamiltonian $H: T^\ast Q\rightarrow
\mathbb{R}$, the external force map $F: T^\ast Q\rightarrow T^\ast
Q$, and the control subset $W$ of \;$T^\ast Q$ are $G$-invariant.
\end{defi}

For a symmetric CH system $(T^\ast Q, G, B, H, F, W)$, if the
$G$-invariant control subset $W$ is a controllability submanifold of
the symmetric CH system, then $W$ is called a {\bf $G$-invariant
controllability submanifold} of this system; if a $G$-invariant
distribution $D\subset TT^*Q|_W$ is a controllability distribution
of the symmetric CH system, then $D$ is called a {\bf $G$-invariant
controllability distribution} of this system, simply written as
$D^G$. In this case $D_W=D^G\cap TW$ is also $G$-invariant
distribution on $W$. For the presheaf $C_{W/D_W}^\infty$, the
$(D^G,D_W)$-local extension property shows that $F$ is a local
$D^G$-invariant extension of $f\circ \pi_{D_W}$ at the point $m\in
\pi_{D_W}^{-1}(V)$, for any open $G$-invariant subset of $W/D_W, $
that is, the local extensions of
 $f\circ \pi_{D_W}$ defined in (2.1) can always be chosen
 to satisfy
\begin{equation}
 \mathrm{d}F(n)|_{D^G(n)}=0,\; \mbox{for any } n\in
 \pi_{D_W}^{-1}(V)\cap U_m.
\end{equation}
From Definition 2.4 it is easy to see that the $G$-invariant
distribution $D^G$ is Poisson if the condition
$\mathrm{d}f|_{D^G}=\mathrm{d}g|_{D^G}=0$, for any $f,g\in
C_{T^*Q}^\infty(U)^G$, $U\subset T^*Q$ an open $G$-invariant subset,
implies that $\mathrm{d}\{f,g\}|_{D^G}=0$. Moreover, assume that the
$G$-invariant distribution $D^G$ is Poisson, and the presheaf
$C_{W/D_W}^\infty$ has the $(D^G,D_W)$-local extension property. If
$F$ and $G\in C_{T^*Q}^\infty(U_m)^G$ are local $D^G$-invariant
extensions of $f\circ \pi_{D_W}$ and $g\circ \pi_{D_W}$ at $m\in
\pi_{D_W}^{-1}(V)$, then their Poisson bracket $\{F,G\}$ is also
local $D^G$-invariant extension.\\

Now for a symmetric CH system $(T^*Q, G, B, H, F, W)$, assume that
$W$ is its $G$-invariant controllability submanifold and $D^G
\subset TT^*Q|_W$ is its $G$-invariant controllability distribution,
such that $D_W=D^G\cap TW$ is a smooth regular $G$-invariant
integrable distribution on $W$, and the presheaf $C_{W/D_W}^\infty$
has the $(D^G,D_W)$-local extension property, then we may consider
that $(T^*Q,B,D^G,W)$ is Poisson reducible if
$(W/D_W,C_{W/D_W}^\infty,B^{W/D_W})$ is a well-defined presheaf of
Poisson algebras. Thus, by using the Poisson reduction by
$G$-invariant controllability distribution $D^G$ for the phase space
of symmetric CH system, we can define the Poisson reducible
symmetric CH system by $G$-invariant controllability distribution as
follows.
\begin{defi}
A symmetric CH system $(T^*Q, G, B, H, F, W)$ is called to be
Poisson reducible by $G$-invariant controllability distribution
$D^G$, if $(T^*Q,B,D^G,W)$ is Poisson reducible by $G$-invariant
distribution $D^G$ in the sense of Definition 2.5, that is,
$(W/D_W,C_{W/D_W}^\infty,B^{W/D_W})$ is a well-defined presheaf of
Poisson algebras, where the Poisson tensor $B^{W/D_W}$ is defined
by, for any open $G$-invariant
  set $V\subset W/D_W$, the Poisson bracket
  $\{\cdot,\cdot\}_V^{W/D_W}: C_{W/D_W}^\infty(V)\times
  C_{W/D_W}^\infty(V)\to C_{W/D_W}^\infty(V)$ is given by
  \begin{equation*}  \label{gpred}
  \{k,l\}_V^{W/D_W}(\pi_{D_W}(m))=\{K,L\}_B(m)
  \end{equation*}
for any $m\in \pi_{D_W}^{-1}(V)$ and any local $D^G$-invariant
extensions $K,L\in C^\infty_{T^*Q}(U_m)^G$ of $k\circ \pi_{D_W}$ and
$l\circ \pi_{D_W}$ at $m$, respectively, where $\pi_{D_W}: W
\rightarrow W/D_W$ is the projection, and $U_m$ is an open
$G$-invariant neighborhood of $m$, such that $U_m\cap W$ is
$D_W$-invariant.
\end{defi}

By using the above definition and Proposition 2.2, we can obtain the
following theorem to give the Poisson reducible conditions of
symmetric CH system by $G$-invariant controllability distribution.

\begin{theo}
Suppose that $(T^*Q, G, B, H, F, W)$ is a symmetric CH system, $W$
is $G$-invariant controllability submanifold and and $D^G \subset
TT^*Q|_W$ is $G$-invariant controllability distribution, such that
$D_W=D^G\cap TW$ is a smooth regular $G$-invariant integrable
distribution on $W$, and the presheaf $C_{W/D_W}^\infty$ has the
$(D^G,D_W)$-local extension property. Then the symmetric CH system
$(T^*Q, G, B, H, F, W)$ is Poisson reducible by $G$-invariant
controllability distribution $D^G$, if and only if for any $z \in W$
  \begin{equation}
    B^\sharp(\Delta^G_z)\subset (\Delta^G_z)_W^\circ\label{prdcon}
  \end{equation}
where
\begin{equation*}\Delta^G_z:=\left\{\mathrm{d}f(z)\left|\begin{aligned}f\in
  C_{T^*Q}^\infty(U_z)^G,\; \mathrm{d}f(y)|_{D^G(y)}=0,\; \mbox{for all } y\in
  U _z\cap W, \mbox{and}\\ \mbox{for any open $G$-invariant neighborhood }U_z \mbox{ of } z \mbox{ in }
  T^*Q\end{aligned}\right.\right\},\end{equation*}
  and \begin{equation*}
    (\Delta^G_z)_W:=\left\{\mathrm{d}f(z)\in \Delta^G_z\left|\begin{aligned}f|_{U_z\cap W}\
    \mbox{is constant for any open $G$-}\\ \mbox{invariant neighborhood } U_z
    \mbox{ of } z \mbox{ in } T^*Q \end{aligned}\right.\right\}.
  \end{equation*}
\end{theo}

\noindent{\bf Proof.} If the symmetric CH system $(T^*Q, G, B, H, F,
W)$ is Poisson reducible by $G$-invariant controllability
distribution $D^G$, then $(W/D_W,C_{W/D_W}^\infty,B_{W/D_W})$ is a
well-defined presheaf of Poisson algebras and for any open
$G$-invariant subset $V\subset W/D_W$, the Poisson tensor
$B^{W/D_W}$ is defined by the following Poisson bracket
\begin{equation}  \label{gpred}
  \{k,l\}_V^{W/D_W}(\pi_{D_W}(z))=\{K,L\}_B(z)
\end{equation}
for any $z\in \pi_{D_W}^{-1}(V)$ and any local $D^G$-invariant
extensions $K,L\in C^\infty_{T^*Q}(U_z)^G$ of $k \cdot \pi_{D_W}$
and $l \cdot \pi_{D_W}$ at $z$, respectively, where $\pi_{D_W}: W
\rightarrow W/D_W$ is the projection, and $U_z$ is an open
$G$-invariant neighborhood of $z$, such that $U_z\cap W$ is
$D_W$-invariant. We shall prove that the inclusion (4.2) holds. For
any $\alpha=\mathrm{d}K(z)\in \Delta^G_{z}, $ that is, $K \in
C_{T^*Q}^\infty(U_{z})^G, $ where $U_{z}$ is an open $G$-invariant
neighborhood of $z$ in $T^*Q, $ such that $U_{z}\cap W$ is
$D_{W}$-invariant and $ \mathrm{d}K(s)|_{D^G(s)}=0,\; \mbox{for all
} s\in U_{z}\cap W. $ In the same way, for any
$\beta=\mathrm{d}L(z)\in (\Delta^G_{z})_{W}, $ that is, $L \in
C_{T^*Q}^\infty(U_{z})^G $ and $L|_{U_{z}\cap W}$ is constant for
the open $G$-invariant neighborhood $U_{z}$ of $z$ in $T^*Q$. Since
$K|_{U_{z}\cap W}$ and $L|_{U_{z}\cap W}$ are $D_W$-invariant, they
have constant values on the leaves of $D_W$. Note that $W$ is
endowed with the relative topology, and $U_{z}\cap W$ is open in
$W$. Thus, we can take that the open set $V:= \pi_{D_W}(U_{z}\cap W)
\subset W/ D_W, $ and define functions $k,\;l \in
C_{W/D_W}^\infty(V)$ by using $K,\;L \in C_{T^*Q}^\infty(U_{z})^G, $
such that $k \cdot \pi_{D_W}|_{\pi_{D_W}^{-1}(V)\cap U_z}=
K|_{\pi_{D_W}^{-1}(V)\cap U_z}$ and $l \cdot
\pi_{D_W}|_{\pi_{D_W}^{-1}(V)\cap U_z}= L|_{\pi_{D_W}^{-1}(V)\cap
U_z}$. Because $L$ is constant on $U_{z}\cap W, $ hence the function
$l$ is constant on the neighborhood $V= \pi_{D_W}(U_{z}\cap W)$ of
$\pi_{D_W}(z)$ in $ W/ D_W$, and we have that
$\{k,l\}_V^{W/D_W}(\pi_{D_W}(z))=0. $ In consequence, from (4.3) we
have that $$<\beta, B^\sharp(\alpha)> = <\mathrm{d}L(z),
B^\sharp(\mathrm{d}K(z))>= \{K, L\}_{B}(z) =
\{k,l\}_V^{W/D_W}(\pi_{D_W}(z))=0. $$ Thus, we prove that
$B^\sharp(\alpha)\in (\Delta^G_{z})_{W}^\circ$ and obtain the
desired inclusion $B^\sharp(\Delta^G_{z})\subset
(\Delta^G_{z})_{W}^\circ .$\\

Conversely, if assume that the inclusion (4.2) holds, we shall prove
that the symmetric CH system $(T^*Q, G, B, H, F, W)$ is Poisson
reducible by $G$-invariant controllability distribution $D^G$, that
is, $(W/D_W,C_{W/D_W}^\infty,B^{W/D_W})$ is a well-defined presheaf
of Poisson algebras. For any open $G$-invariant subset $V\subset
W/D_W$, and functions $k,\;l \in C_{W/D_W}^\infty(V)$, as well as
$K,\;L \in C_{T^*Q}^\infty(U_{z})^G, $ which are the local
$D^G$-invariant extensions of $k \cdot \pi_{D_W}$ and $l \cdot
\pi_{D_W}$ at $z\in \pi_{D_W}^{-1}(V)$, respectively. We define the
Poisson tensor $B^{W/D_W}$ by the Poisson bracket
$\{\cdot,\cdot\}_V^{W/D_W}: C_{W/D_W}^\infty(V)\times
C_{W/D_W}^\infty(V)\to C_{W/D_W}^\infty(V), $ which is given by
(4.3). We shall prove that this Poisson bracket has the property of
local $D^G$-invariant extension, that is,
$$\{k,l\}_V^{W/D_W} \cdot \pi_{D_W}|_{\pi_{D_W}^{-1}(V)\cap U_z}=
\{K, L\}_{B}|_{\pi_{D_W}^{-1}(V)\cap U_z}. $$ In order to do this,
we have to check that $\{K, L\}_{B}(z)$ doesn't depend on the choice
of the point $z\in \pi_{D_W}^{-1}(V)$ and the local extensions $K$
and $L$. Because $D^G$ is a $G$-invariant controllability
distribution, it is a $G$-invariant Poisson integrable generalized
distribution, if the presheaf $C_{W/D_W}^\infty$ has the
$(D^G,D_W)$-local extension property, and $K$ and $L\in
C_{T^*Q}^\infty(U_z)^G$ are local $D^G$-invariant extensions of
$k\circ \pi_{D_W}$ and $l\circ \pi_{D_W}$ at $z\in
\pi_{D_W}^{-1}(V)$, then their Poisson bracket $\{K, L\}$ is also
local $D^G$-invariant extension. Thus the function $\{K, L\}|_{U_z
\cap W}$ is constant along the integral curves of any section of
$D_W$. Take $z,\; z' \in T^*Q, $ such that $\pi_{D_W}(z)=
\pi_{D_W}(z'). $ Assume that $K',\;L' \in
C_{T^*Q}^\infty(U_{z'})^G,$ are the local $D^G$-invariant extensions
of $k \cdot \pi_{D_W}$ and $l \cdot \pi_{D_W}$ at $z'\in
\pi_{D_W}^{-1}(V)$, respectively, where $U_{z'}$ is an open
$G$-invariant neighborhood of $z'$ in $T^*Q, $ such that $U_{z'}\cap
W$ is $D_{W}$-invariant. Since $z'$ can be connected to $z$ by a
finite union of integral curves of sections of $D_W,$ and $U_{z}\cap
W$ and $U_{z'}\cap W$ are both $D_{W}$-invariant, it follows that
$U_{z}\cap U_{z'}\cap W$ contains $z$ and $z'$ and it is also
$D_{W}$-invariant. Thus, $\{K', L'\}_{B}(z)= \{K', L'\}_{B}(z'). $
Next, we shall check that $\{K', L'\}_{B}(z)= \{K, L\}_{B}(z), $
that is, that $\{K, L\}_{B}(z)$ doesn't depend on the choice of the
extensions $K$ and $L$. Because of the antisymmetry of $\{\cdot,
\cdot\}_{B}, $ it suffices to show that it doesn't depend on the
choice of the extension $L. $ Let $L' \in C_{T^*Q}^\infty(U_{z})^G,
$ be another local $D^G$-invariant extension of $l \cdot \pi_{D_W}$.
Then we have that $(L-L')|_{\pi_{D_W}^{-1}(V)\cap U_z} =0, $ and
$\mathrm{d}(L-L')|_{D^G(s)}= \mathrm{d}L|_{D^G(s)}-
\mathrm{d}L'|_{D^G(s)}=0, $ for any $s\in U_{z}\cap W. $ Hence
$\mathrm{d}(L-L')(z)\in (\Delta^G_{z})_{W}. $ Notice that
$\mathrm{d}K(z)\in \Delta^G_{z}$ by definition, and by using the
inclusion (4.2), we have that $$\{K, L- L'\}_{B}(z)=
<\mathrm{d}(L-L')(z), B^\sharp(\mathrm{d}K(z))>= 0, $$ that is,
$\{K, L\}_{B}(z)= \{K, L'\}_{B}(z). $ At last, the Leibniz and
Jacobi identities for $\{\cdot, \cdot\}_V^{W/D_W}$ follow directly
from the definition of $\{\cdot, \cdot\}_V^{W/D_W}$ and the fact
that $\{\cdot, \cdot\}_{B}$ satisfies these identities. $\hskip 1cm
\blacksquare$

\section{Singular Poisson Reduction of CH System with Symmetry }

In this section, we shall introduce the singular Poisson reduction
of CH system with symmetry $(T^*Q, G, B, H, F, W)$ and its
$SP$-reduced system $(M^{(K)}, B^{(K)}, h^{(K)}, f^{(K)}, W^{(K)})$,
as well as SPR-CH-equivalence by using the controlled Hamiltonian
method given in Marsden et al\cite{mawazh10}. These are a
generalization of the regular Poisson reduction of CH system with
symmetry given in \cite {wazh12} to the singular context. we shall
also follow the notations and conventions for singular reduction of
a differential space introduced in Jotz, Ratiu and
$\acute{S}$niatycki \cite{jorasn11}, Pflaum \cite{pf01},
Sjamaar and Lerman \cite{sjle91}.\\

Let $Q$ be a smooth manifold and $T^\ast Q$ its cotangent bundle
with a associated Poisson bracket $\{\cdot, \cdot \}$ and Poisson
tensor $B$. Let $\Phi:G\times Q\rightarrow Q$ be a smooth left
action of the Lie group $G$ on $Q$, and the cotangent lifted left
action $\Phi^{T^\ast}:G\times T^\ast Q\rightarrow T^\ast Q$ is
canonical, and proper, but may not be free. Thus, the orbit space
$T^\ast Q /G$ is not necessarily smooth manifold, but just a
stratified topological space, and the projection $\pi_{/G}: T^\ast Q
\rightarrow T^\ast Q /G$ is a surjective submersion. In the
following we shall describe the structure of the orbit space $T^\ast
Q /G$. For a closed Lie subgroup $K$ of $G$, we define the isotropy
type set $(T^*Q)_K = \{m\in T^*Q | G_m =K \}, $ where $G_m= \{g\in G
| gm=m \} $ is the isotropy subgroup of $m\in T^*Q. $ Since the
$G$-action on $T^*Q$ is proper, all isotropy groups are compact, and
the sets $(T^*Q)_K$, where $K$ ranges over the closed Lie subgroups
of $G$ for which $(T^*Q)_K$ is nonempty, form a partition of $T^*Q$.
Define the normalizer of $K$ in $G$ by $N(K):= \{g\in G | gKg_{-1}=K
\}, $ $N(K)$ is a closed Lie subgroup of $G$. Since $K$ is a normal
subgroup of $N(K)$, the quotient group $N(K)/ K$ is a Lie group. If
$m\in (T^*Q)_K$, we have that $G_m =K, $ and for all $g\in G, \;
G_{gm}= gKg^{-1}. $ Thus, $gm$ lies in $(T^*Q)_K $ if and only if
$g\in N(K), $ and the action of $G$ on $T^*Q$ restricts to an action
of $N(K)$ on $(T^*Q)_K, $ which induces a free and proper action of
$N(K)/ K$ on $(T^*Q)_K. $ Define the orbit type set $(T^*Q)_{(K)}=
\{m\in T^*Q | G_m \in (K) \}, $ where $(K)$ is set of $K$ conjugate
classes. Then, $(T^*Q)_{(K)}= \{gm | g\in G, \; m\in (T^*Q)_{K}\} =
\pi_{/G}^{-1}(\pi_{/G}((T^*Q)_{K})). $ In these cases, the connected
components of $(T^*Q)_{K}$ and $(T^*Q)_{(K)}$ are embedded
submanifolds of $T^*Q, $ therefore $(T^*Q)_{K}$ is an isotropy type
submanifold and $(T^*Q)_{(K)}$ is an orbit type submanifold.
Moreover, $\pi_{/G}((T^*Q)_{(K)})= \{gm | m\in (T^*Q)_{K}\}/G =
(T^*Q)_{K}/ N(K)= (T^*Q)_{K}/ (N(K)/ K). $ But the action of $N(K)/
K$ on $(T^*Q)_{K}$ is free and proper, which implies that
$(T^*Q)_{K}/ (N(K)/ K)$ is a quotient manifold. Thus,
$\pi_{/G}((T^*Q)_{(K)})$ is a manifold contained in the orbit space
$T^*Q /G. $ From Jotz et al \cite{jorasn11} we know that a partition
of the orbit space $T^*Q /G$ by connected components of
$\pi_{/G}((T^*Q)_{(K)})$ is a decomposition of $T^*Q /G$ as a
differential space. The corresponding stratification of $T^*Q /G$ is
called the orbit type stratification of the orbit space, and the
strata, denoted by $M^{(K)}$, which are the connected components of
$\pi_{/G}((T^*Q)_{(K)})$ of such stratification, are orbits of the
family of all vector fields on $T^*Q /G$. Moreover, the each stratum
$M^{(K)}$ is a Poisson manifold with the Poisson bracket $\{\cdot,
\cdot \}^{(K)}$ uniquely characterized by the relation
\begin{equation}\{f_K, g_K \}^{(K)}\cdot (\pi^{(K)}(\alpha))= \{f_K\cdot \pi^{(K)},
g_K\cdot \pi^{(K)}\}_B(\alpha),\;\; \forall \alpha \in T^\ast Q,\;\;
f_K, g_K \in C^{\infty}(M^{(K)}),
\end{equation}
which is a induced Poisson structure by that of $T^*Q$,
such that the projection $\pi^{(K)}: T^\ast Q \rightarrow M^{(K)} $
is a Poisson map. On the other hand, from Ortega and Ratiu
\cite{orra04}, we know that for the free and proper $G$-action on
$T^\ast Q$, the orbit space $(T^\ast Q)/G$ is Poisson diffeomorphic
to a Poisson fiber bundle with respect to its natural induced
Poisson bracket, if a connection $\mathcal{A}$ on the $G$-principal
bundle $\pi: Q \rightarrow Q/G $ is introduced. We assume that the
bundle structures of $T^*Q$ and $(T^\ast Q)/G$ is compatible with
the stratification of $T^*Q /G$, such that the each stratum
$M^{(K)}$ is a Poisson fiber bundle with Poisson tensor $B^{(K)}$,
where $B^{(K)}$ denotes the Poisson tensor induced by the bracket
$\{\cdot, \cdot \}^{(K)}$ on $M^{(K)}$. From (5.1) we have that
$\pi^{(K)^\ast} B^{(K)}= B$. The pair $(M^{(K)}, B^{(K)})$ is called
the {\bf singular Poisson reduced stratum} of $(T^\ast Q, B)$.\\

For the Hamiltonian systems associated to $G$-invariant and
stratified Hamiltonian functions, we can naturally reduce to
Hamiltonian systems on these strata. In fact, If $H: T^\ast Q
\rightarrow \mathbb{R}$ is a $G$-invariant and stratified
Hamiltonian, the flow $F_t$ of the Hamiltonian vector field $X_H$
leaves the orbit type stratification of the orbit space $T^*Q /G $
invariant and commutes with the $G$-action, so it induces a flow
$f_t^{(K)}$ on each stratum $M^{(K)}$, defined by $f_t^{(K)}\cdot
\pi^{(K)}=\pi^{(K)} \cdot F_t $, and the vector field $X_{h^{(K)}}$
generated by the flow $f_t^{(K)}$ on $(M^{(K)}, B^{(K)})$ is
Hamiltonian with the associated singular Poisson reduced Hamiltonian
function $h^{(K)}: M^{(K)} \rightarrow \mathbb{R}$ defined by
$h^{(K)}\cdot \pi^{(K)}= H $ and the Hamiltonian vector fields $X_H$
and $X_{h^{(K)}}$ are $\pi^{(K)}$-related. Moreover, if consider
that the external force map $F: T^\ast Q \rightarrow T^\ast Q$ and
the controlled subset $W$ are $G$-invariant and stratified, then we
can introduce a kind of the singular Poisson reducible CH systems as
follows.
\begin{defi}
(Singular Poisson Reducible CH system) A 6-tuple $(T^\ast Q, G, B,
H, F, W)$, where the Hamiltonian $H: T^\ast Q\rightarrow
\mathbb{R}$, the fiber-preserving map $F: T^\ast Q \rightarrow
T^\ast Q$ and the fiber submanifold $W$ of $T^\ast Q$ are all
$G$-invariant and stratified, is called a {\bf singular Poisson
reducible CH system}, if for isotropy subgroup $K (\subset G) $ of
the $G$-action on $T^*Q$, the connected components of
$\pi_{/G}((T^*Q)_{(K)})$ form a corresponding orbit type
stratification of the orbit space $T^*Q /G$, and there exists the
singular Poisson reduced system, that is, the 5-tuple $(M^{(K)},
B^{(K)}, h^{(K)}, f^{(K)}, W^{(K)})$, where the stratum $M^{(K)}$ is
a connected component of the projection $\pi_{/G}((T^*Q)_{(K)})$ of
the $(K)$-orbit type submanifold $(T^\ast Q)_{(K)}$ to the orbit
space, $\pi^{(K)\ast} B^{(K)}= B$,\; $h^{(K)}\cdot \pi^{(K)}= H $,\;
$f^{(K)}\cdot \pi^{(K)}=\pi^{(K)}\cdot F$,\; $W^{(K)}=\pi^{(K)}(W)$,
is a CH system, which is simply written as $SP$-reduced CH system.
Where $(M^{(K)}, B^{(K)})$ is the singular Poisson reduced stratum,
the function $h^{(K)}: M^{(K)} \rightarrow \mathbb{R}$ is called the
reduced Hamiltonian, the fiber-preserving map $f^{(K)}: M^{(K)}
\rightarrow M^{(K)} $ is called the reduced (external) force map,
$W^{(K)}$ is a fiber submanifold of $M^{(K)}$ and is called the
reduced control subset.
\end{defi}

Denote by $X_{(T^\ast Q,G,B,H,F,u)}$ the vector field of singular
Poisson reducible CH system $(T^\ast Q,G,B, \\H,F,W)$ with a control
law $u$, then
\begin{equation}X_{(T^\ast Q,G,B,H,F,u)}=
B^\sharp\mathrm{d}H+\textnormal{vlift}(F)+\textnormal{vlift}(u)
\end{equation}
Moreover, for the singular Poisson reducible CH system we can also
introduce the singular Poisson reduced controlled Hamiltonian
equivalence (SPR-CH-equivalence) as follows.

\begin{defi}(SPR-CH-equivalence)
Suppose that we have two singular Poisson reducible CH systems
$(T^\ast Q_i, G_i,B_i,H_i, F_i, W_i),\; i=1,2,$ we say them to be
{\bf SPR-CH-equivalent} or simply \\$(T^\ast Q_1,
G_1,B_1,H_1,F_1,W_1)\stackrel{SPR-CH}{\sim}(T^\ast
Q_2,G_2,B_2,H_2,F_2,W_2)$, if there exists a diffeomorphism
$\varphi:Q_1\rightarrow Q_2$, such that the following Hamiltonian
matching conditions hold:

\noindent {\bf SPR-H1:} The cotangent lift map $\varphi^\ast: T^\ast
Q_2 \rightarrow T^\ast Q_1$ is a Poisson map, and $W_1=\varphi^\ast
(W_2)$;

\noindent {\bf SPR-H2:} For the Lie group actions $\Phi_i^{T^\ast}:
G_i\times T^\ast Q_i\rightarrow T^\ast Q_i,\; i=1,2,$ the map
$\varphi^\ast: T^\ast Q_2\rightarrow T^\ast Q_1$ is
$(G_2,G_1)$-equivariant, and for the corresponding orbit type
stratification of the orbit spaces $T^*Q_i /G_i,\; i=1,2, $ the
induced orbit space map $\varphi_{/G}^{\ast}: T^*Q_2 /G_2
\rightarrow T^*Q_1 /G_1$ is stratified;

\noindent {\bf SPR-H3:}
$Im[B_1^\sharp\mathrm{d}H_1+\textnormal{vlift}(F_1)-T \varphi^\ast
(B_2^\sharp \mathrm{d}H_2)-\textnormal{vlift}(\varphi^\ast
F_2\varphi_\ast)]\subset\textnormal{vlift}(W_1)$.
\end{defi}

It is worthy of note that for the singular Poisson reducible CH
system, the equivalent map not only keeps the Poisson structure, but
also keeps the equivariance of G-action on the cotangent bundle and
the stratification. If a feedback control law $u^{(K)}: M^{(K)}
\rightarrow W^{(K)}$ is chosen, the $SP$-reduced CH system
$(M^{(K)}, B^{(K)}, h^{(K)}, f^{(K)}, W^{(K)})$ is a closed-loop
dynamic system with a control law $u^{(K)}$. Assume that its vector
field $X_{(M^{(K)}, B^{(K)}, h^{(K)}, f^{(K)}, u^{(K)})}$ can be
expressed by
\begin{equation} X_{(M^{(K)}, B^{(K)}, h^{(K)},
f^{(K)}, u^{(K)})}= (B^{(K)})^{\sharp} \mathrm{d}h^{(K)}+
\textnormal{vlift}(f^{(K)})+\textnormal{vlift}(u^{(K)}),
\end{equation}
where $(B^{(K)})^{\sharp} \mathrm{d}h^{(K)}= X_{h^{(K)}}$,
$\textnormal{vlift}(f^{(K)})=
\textnormal{vlift}(f^{(K)})X_{h^{(K)}}$,
$\textnormal{vlift}(u^{(K)})=
\textnormal{vlift}(u^{(K)})X_{h^{(K)}}$, and satisfies the condition
that $X_{(M^{(K)}, B^{(K)}, h^{(K)}, f^{(K)}, u^{(K)})}$ and
$X_{(T^\ast Q,G,B,H,F,u)}$ are $\pi^{(K)}$-related, that is,
\begin{equation}X_{(M^{(K)}, B^{(K)}, h^{(K)},
f^{(K)}, u^{(K)})}\cdot \pi^{(K)}=T\pi^{(K)}\cdot X_{(T^\ast
Q,G,B,H,F,u)},
\end{equation}
where $T\pi^{(K)}: TT^*Q \rightarrow TM^{(K)}$ is tangent map of the
projection $\pi^{(K)}: T^*Q \rightarrow M^{(K)}. $ Then we can
obtain the following singular Poisson reduction theorem for CH
system, which explains the relationship between the
SPR-CH-equivalence for the singular Poisson reducible CH systems
with symmetry and the CH-equivalence for associated $SP$-reduced CH
systems.
\begin{theo} Two singular Poisson reducible CH systems
$(T^\ast Q_i, G_i, B_i, H_i, F_i, W_i)$, $i=1,2,$ are
SPR-CH-equivalent if and only if the
associated $SP$-reduced CH systems $(M_i^{(K_i)}, B_i^{(K_i)}, h_i^{(K_i)},\\
f_i^{(K_i)}, W_i^{(K_i)}), i=1,2$, are CH-equivalent.
\end{theo}

\noindent{\bf Proof.} If $(T^\ast Q_1, G_1, B_1, H_1, F_1, W_1)
\stackrel{SPR-CH}{\sim}(T^\ast Q_2, G_2, B_2, H_2, F_2, W_2)$, then
there exists a diffeomorphism $\varphi: Q_1\rightarrow Q_2$ such
that $\varphi^\ast: T^\ast Q_2\rightarrow T^\ast Q_1$ is a
$(G_2,G_1)$-equivariant Poisson map, $W_1=\varphi^\ast (W_2)$, and
for the corresponding orbit type stratification of the orbit spaces
$T^*Q_i /G_i,\; i=1,2, $ the induced orbit space map
$\varphi_{/G}^{\ast}: T^*Q_2 /G_2 \rightarrow T^*Q_1 /G_1$ is
stratified, and {\bf SPR-H3} holds. From the following commutative
Diagram-11:
$$\;\;\;\;\;\;
\begin{CD}
T^\ast Q_2/G_2 @<\pi_{/G_2}<< T^\ast Q_2  @>\pi^{(K_2)}>> M_{2}^{(K_2)}\\
@V\varphi^\ast_{/G} VV  @V\varphi^\ast VV  @V\varphi^\ast_{(K)} VV \\
T^\ast Q_1/G_1 @<\pi_{/G_1}<< T^\ast Q_1  @>\pi^{(K_1)}>>
M_{1}^{(K_1)}
\end{CD}
\;\;\;\;\;\;\;\;\;\;\;\;\;\;\;\;\;\;\;\;\;
\begin{CD}
\Lambda^2(T^*T^\ast Q_1) @<\pi^{(K_1)\ast} << \Lambda^2(T^*M_{1}^{(K_1)})\\
 @V(\varphi^\ast)^\ast VV  @V(\varphi^\ast_{(K)})^{\ast} VV \\
\Lambda^2(T^*T^\ast Q_2) @<\pi^{(K_2)\ast} <<
\Lambda^2(T^*M_{2}^{(K_2)})
\end{CD}
$$
$$\;\;\;\;\;\mbox{Diagram-11}\;\;\;\;\;\;\;\;\;\;\;\;\;\;\;\;\;\;\;\;\;\;\;\;\;\;\;
\;\;\;\;\;\;\;\;\;\;\;\;\;\;\;\;\;\;\;\;\;\;\;\;\;
\mbox{Diagram-12}$$ We can define a map $\varphi_{(K)}^{\ast}:
M_{2}^{(K_2)} \rightarrow M_{1}^{(K_1)} $ such that
$\varphi_{(K)}^{\ast} \cdot \pi^{(K_2)}=\pi^{(K_1)}\cdot
\varphi^\ast $, and $\varphi_{/G}^\ast \cdot
\pi_{/G_2}=\pi_{/G_1}\cdot \varphi^\ast $. Because $\varphi^\ast:
T^\ast Q_2\rightarrow T^\ast Q_1$ is $(G_2,G_1)$-equivariant and the
induced orbit space map $\varphi^\ast_{/G}: T^\ast Q_2 /G_2
\rightarrow T^\ast Q_1/G_1$ is stratified, the map
$\varphi_{(K)}^{\ast} = \varphi^\ast_{/G} |_{M^{(K)}}: M_{2}^{(K_2)}
\rightarrow M_{1}^{(K_1)}$ is well-defined. We shall show that
$\varphi_{(K)}^\ast$ is a Poisson map and
$W_{1}^{(K_1)}=\varphi_{(K)}^\ast (W_{2}^{(K_2)})$. In fact, since
$\varphi^\ast: T^\ast Q_2 \rightarrow T^\ast Q_1$ is a Poisson map,
the map $(\varphi^\ast)^\ast: \Lambda^2 (T^*T^\ast Q_1)\rightarrow
\Lambda^2 (T^*T^\ast Q_2)$ such that $(\varphi^\ast)^\ast B_1= B_2$,
and by using (5.1), the map $\pi^{(K_i)}: T^\ast Q_i \rightarrow
M_{i}^{(K_i)}$ is also a Poisson map, $B_i=\pi^{(K_i)\ast}
B_{i}^{(K_i)},\; i=1,2$, from the commutative Diagram-12, we have
that
\begin{align*}
\pi^{(K_2)\ast}\cdot(\varphi_{(K)}^\ast)^\ast B_{1}^{(K_1)}&=(
\varphi_{(K)}^{\ast} \cdot \pi^{(K_2)})^\ast
B_{1}^{(K_1)}=(\pi^{(K_1)}\cdot \varphi^\ast)^\ast
B_{1}^{(K_1)}\\&=(\varphi^\ast)^\ast \cdot \pi^{(K_1)\ast}
B_{1}^{(K_1)} =(\varphi^\ast)^\ast B_1= B_2 =\pi^{(K_2)\ast}\cdot
B_{2}^{(K_2)}.
\end{align*}
Notice that $\pi^{(K_2)\ast}$ is surjective, thus,
$(\varphi_{(K)}^\ast)^\ast B_{1}^{(K_1)}= B_{2}^{(K_2)}$. Because $
W_{i}^{(K_i)}=\pi^{(K_i)}(W_i),\; i=1,2$ and $W_1=\varphi^\ast
(W_2)$, we have that
$$W_{1}^{(K_1)}=\pi^{(K_1)}(W_1)=\pi^{(K_1)}\cdot\varphi^\ast
(W_2)=\varphi_{(K)}^\ast\cdot
\pi^{(K_2)}(W_2)=\varphi_{(K)}^\ast(W_{2}^{(K_2)}).$$ Next, from
(5.2) and (5.3), we know that for $i=1,2$,
$$X_{(T^\ast Q_i, G_i, B_i, H_i, F_i,
u_i)}=B_i^\sharp\mathrm{d}H_i+\textnormal{vlift}(F_i)+\textnormal{vlift}(u_i),$$
$$X_{(M_{i}^{(K_i)}, B_{i}^{(K_i)}, h_{i}^{(K_i)},
f_{i}^{(K_i)},u_{i}^{(K_i)})}=(B_{i}^{(K_i)})^\sharp\mathrm{d}h_{i}^{(K_i)}
+\textnormal{vlift}(f_{i}^{(K_i)})+\textnormal{vlift}(u_{i}^{(K_i)}),$$
and from (5.4), we have that
$$X_{(M_{i}^{(K_i)}, B_{i}^{(K_i)}, h_{i}^{(K_i)},
f_{i}^{(K_i)},u_{i}^{(K_i)})}\cdot \pi^{(K_i)}=T\pi^{(K_i)}\cdot
X_{(T^\ast Q_i, G_i, B_i, H_i, F_i, u_i)}.$$ Since $H_i,F_i$ and
$W_i$ are all $G_i$-invariant and stratified, $i=1,2$, and
$$h_{i}^{(K_i)}\cdot \pi^{(K_i)}=H_i ,\;\;f_{i}^{(K_i)}\cdot \pi^{(K_i)}=\pi^{(K_i)}\cdot F_i,\;\;
u_{i}^{(K_i)}\cdot \pi^{(K_i)}= \pi^{(K_i)}\cdot u_i, \;\;i=1,2.$$
From the following commutative Diagram-13,
\[
\begin{CD}
T^\ast Q_2 @>\textnormal{vlift}>> TT^\ast Q_2 @>T\pi^{(K_2)}
>> T(M_{2}^{(K_2)}) @<\textnormal{vlift} << M_{2}^{(K_2)}\\
@V\varphi^\ast VV @V T\varphi^\ast VV @V T\varphi^\ast_{(K)} VV @V\varphi_{(K)}^\ast VV \\
T^\ast Q_1 @>\textnormal{vlift}>> TT^\ast Q_1 @>T\pi^{(K_1)} >>
T(M_{1}^{(K_1)}) @<\textnormal{vlift} << M_{1}^{(K_1)}
\end{CD}
\]
$$\mbox{Diagram-13}$$
we have that
$$\textnormal{vlift}(\varphi_{(K)}^\ast\cdot
f_{2}^{(K_2)}\cdot \varphi_{(K) \ast})\cdot \pi^{(K_1)}=
T\pi^{(K_1)}\cdot \textnormal{vlift}(\varphi^\ast
F_2\varphi_\ast),$$
$$\textnormal{vlift}(\varphi_{(K)}^\ast\cdot
u_{2}^{(K_2)}\cdot \varphi_{(K) \ast})\cdot \pi^{(K_1)}=
T\pi^{(K_1)}\cdot \textnormal{vlift}(\varphi^\ast
u_2\varphi_\ast),$$
$$T\varphi_{(K)}^\ast\cdot ((B_{2}^{(K_2)})^\sharp\mathrm{d}h_{2}^{(K_2)})\cdot \pi^{(K_1)}
=T\pi^{(K_1)}\cdot T\varphi^\ast \cdot(B_2^\sharp \mathrm{d}H_2),$$
where $\varphi_{(K)\ast}=(\varphi^\ast_{(K)})^{-1}: M_{1}^{(K_1)}
\rightarrow M_{2}^{(K_2)}$. From Hamiltonian matching condition {\bf
SPR-H3}, we have that
\begin{equation}Im[(B_{1}^{(K_1)})^\sharp\mathrm{d}h_{1}^{(K_1)}+\textnormal{vlift}(f_{1}^{(K_1)})-
T\varphi_{(K)}^\ast
((B_{2}^{(K_2)})^\sharp\mathrm{d}h_{2}^{(K_2)})-\textnormal{vlift}(\varphi_{(K)}^\ast
\cdot f_{2}^{(K_2)}\cdot {\varphi_{(K)}}_\ast)]\\ \subset
\textnormal{vlift}(W_{1}^{(K_1)}). \end{equation} So,
$$(M_{1}^{(K_1)}, B_{1}^{(K_1)}, h_{1}^{(K_1)}, f_{1}^{(K_1)}, W_{1}^{(K_1)})
\stackrel{CH}{\sim}( M_{2}^{(K_2)}, B_{2}^{(K_2)}, h_{2}^{(K_2)},
f_{2}^{(K_2)}, W_{2}^{(K_2)}).$$

Conversely, assume that the $SP$-reduced CH systems $(M_{i}^{(K_i)},
B_{i}^{(K_i)}, h_{i}^{(K_i)}, f_{i}^{(K_i)}, W_{i}^{(K_i)}), \;
i=1,2,$ are CH-equivalent. Then there exists a diffeomorphism
$\varphi_{(K)}^\ast: M_{2}^{(K_2)} \rightarrow M_{1}^{(K_1)} $,
which is a Poisson map, $W_{1}^{(K_1)}=\varphi_{(K)}^\ast
(W_{2}^{(K_2)})$, and $(5.5)$ holds. We can define the maps
$\varphi^\ast: T^\ast Q_2\rightarrow T^\ast Q_1$, and
$\varphi^\ast_{/G}: T^\ast Q_2 /G_2 \rightarrow T^\ast Q_1/G_1$,
such that $\pi^{(K_1)}\cdot \varphi^\ast= \varphi_{(K)}^\ast\cdot
\pi^{(K_2)}$, and $\varphi_{/G}^\ast \cdot
\pi_{/G_2}=\pi_{/G_1}\cdot \varphi^\ast $, and $ \varphi^\ast_{/G}
|_{M^{(K)}}= \varphi_{(K)}^{\ast}: M_{2}^{(K_2)} \rightarrow
M_{1}^{(K_1)} $, see the commutative Diagram-11, as well as a
diffeomorphism $\varphi:Q_1\rightarrow Q_2$ whose cotangent lift is
just $\varphi^\ast: T^\ast Q_2 \rightarrow T^\ast Q_1$. From
definition of $\varphi^\ast$, we know that $\varphi^\ast$ is
$(G_2,G_1)$-equivariant and the induced orbit space map
$\varphi^\ast_{/G}: T^\ast Q_2 /G_2 \rightarrow T^\ast Q_1/G_1$ is
stratified. In fact, for any $\alpha_i\in T^\ast Q_i$, $g_i\in G_i$,
$i=1,2,$ such that $\alpha_1=\varphi^\ast(\alpha_2)$,
$[\alpha_1]=\varphi^\ast_{(K)}[\alpha_2]$, then we have that
\begin{align*}
\pi^{(K_1)}\cdot\varphi^\ast\cdot(\Phi_{2g_2}(\alpha_2))
&=\pi^{(K_1)}\cdot\varphi^\ast(g_2\alpha_2)
=\varphi_{(K)}^\ast\cdot \pi^{(K_2)}(g_2\alpha_2)=\varphi_{(K)}^\ast[\alpha_2]=[\alpha_1]\\
&=\pi^{(K_1)}(g_1\alpha_1)=\pi^{(K_1)}\cdot(\Phi_{1g_1}(\alpha_1))
=\pi^{(K_1)}\cdot \Phi_{1g_1}\cdot \varphi^\ast(\alpha_2).
\end{align*}
Since $\pi^{(K_1)}$ is surjective, so,
$\varphi^\ast\cdot\Phi_{2g_2}=\Phi_{1g_1}\cdot \varphi^\ast, $ and
hence $\varphi^\ast$ is $(G_2,G_1)$-equivariant. Notice that
$M_{i}^{(K_i)}$ form the corresponding orbit type stratification of
the orbit space $T^*Q_i /G_i,\; i=1,2, $ and by construction $
\varphi^\ast_{/G} |_{M^{(K)}}= \varphi_{(K)}^{\ast}: M_{2}^{(K_2)}
\rightarrow M_{1}^{(K_1)} $. Thus, in this case the induced orbit
space map $\varphi^\ast_{/G}: T^\ast Q_2 /G_2 \rightarrow T^\ast
Q_1/G_1$ is stratified. Moreover,
$\pi^{(K_1)}(W_1)=W_{1}^{(K_1)}=\varphi_{(K)}^\ast(W_{2}^{(K_2)})=\varphi_{(K)}^\ast\cdot
\pi^{(K_2)}(W_2) =\pi^{(K_1)}\cdot \varphi^\ast (W_2)$. Since
$\pi^{(K_1)}$ is surjective, then $W_1=\varphi^\ast (W_2)$. We shall
show that $\varphi^\ast$ is a Poisson map. Because
$\varphi_{(K)}^\ast: M_{2}^{(K_2)} \rightarrow M_{1}^{(K_1)} $ is a
Poisson map, the map $(\varphi_{(K)}^\ast)^\ast:
\Lambda^2(T^*M_{1}^{(K_1)})\rightarrow \Lambda^2(T^*M_{2}^{(K_2)})$
such that $(\varphi_{(K)}^\ast)^\ast B_{1}^{(K_1)}= B_{2}^{(K_2)}$
and by using $(5.1)$, $\pi^{(K_i)}: T^\ast Q_i \rightarrow
M_{i}^{(K_i)}$ is also a Poisson map, $B_i=\pi^{(K_i)\ast}
B_{i}^{(K_i)},\; i=1,2$, from the commutative Diagram-12, we have
that
\begin{align*}
B_2=\pi^{(K_2)\ast}
B_{2}^{(K_2)}=\pi^{(K_2)\ast}\cdot(\varphi_{(K)}^\ast)^\ast
B_{1}^{(K_1)} =(\varphi^\ast)^\ast\cdot \pi^{(K_1)\ast}
B_{1}^{(K_1)} =(\varphi^\ast)^\ast B_1.
\end{align*}
Since the vector field $X_{(T^\ast Q_i,G_i,B_i,H_i,F_i,u_i)}$ and
$X_{(M_{i}^{(K_i)},
B_{i}^{(K_i)},h_{i}^{(K_i)},f_{i}^{(K_i)},u_{i}^{(K_i)})}$ is
$\pi^{(K_i)}$-related, $i=1,2,$ and $H_i, F_i$ and $W_i$ are all
$G_i$-invariant and stratified, $i=1,2,$ in the same way, from
$(5.5)$, we have that
$$Im[B_1^\sharp\mathrm{d}H_1+\textnormal{vlift}(F_1)-T \varphi^\ast
(B_2^\sharp\mathrm{d}H_2)-\textnormal{vlift}(\varphi^\ast
F_2\varphi_\ast)]\subset \textnormal{vlift}(W_1), $$ that is,
Hamiltonian matching condition {\bf SPR-H3} holds. Thus,
$$(T^\ast Q_1, G_1, B_1, H_1, F_1,
W_1)\stackrel{SPR-CH}{\sim}(T^\ast Q_2, G_2, B_2, H_2, F_2,
W_2).\hskip 1cm \blacksquare$$

\begin{rema}
When the $G$-action on $T^*Q$ is free and proper, then $T^*Q /G$ is
a smooth manifold. Itself is only one stratified space. In this
case, from Definition 5.1 and Definition 5.2 we can obtain the
regular reducible CH system and RPR-CH-equivalence, and from Theorem
5.3 we can get the regular Poisson reduction theorem. Thus, we can
recover the regular Poisson reduction results of CH systems given in
Wang and Zhang \cite {wazh12}.
\end{rema}

\section{Relationship Between Singular and by Controllability Distribution Poisson Reductions}

For a singular Poisson reducible CH system $(T^*Q, G, B, H, F, W)$,
it is a symmetric CH system and $W$ is a $G$-invariant submanifold
of $T^*Q$. Assume that $W$ is a controllability submanifold and
$D^G\subset TT^*Q|_W $ is a $G$-invariant controllability
distribution, such that $D_W=D^G\cap TW$ is a smooth regular
$G$-invariant integrable distribution on $W$, and the presheaf
$C_{W/D_W}^\infty$ has the $(D^G,D_W)$-local extension property,
then singular Poisson reducible CH system $(T^*Q, G, B, H, F, W)$ is
also Poisson reducible by $G$-invariant controllability distribution
$D^G$. On the other hand, for the $SP$-reduced CH system $(M^{(K)},
B^{(K)}, h^{(K)}, f^{(K)}, W^{(K)})$, if $W^{(K)}= \pi^{(K)}(W)$ is
a controllability submanifold and $D^{(K)}= T\pi^{(K)}(D^G) \subset
TM^{(K)}|_{W^{(K)}}$ is controllability distribution, then we may
consider that the $SP$-reduced CH system $(M^{(K)}, B^{(K)},
h^{(K)}, f^{(K)}, \\ W^{(K)})$ is Poisson reducible by the reduced
controllability distribution $D^{(K)}$. Moreover, we can also prove
that the property of Poisson reduction for singular Poisson
reducible CH systems by $G$-invariant controllability distribution
leaves invariant under SPR-CH-equivalence conditions.\\

At first, we give the following theorem to state the relationship
between Poisson reduction for singular Poisson reducible CH systems
by $G$-invariant controllability distribution $D^G$ and Poisson
reduction for assocated $SP$-reduced CH system by controllability
distribution $D^G/ G$.

\begin{theo}
Suppose that symmetric CH system $(T^\ast Q,G,B,H,F,W)$ is a
singular Poisson reducible CH system and its associated $SP$-reduced
CH system is $(M^{(K)}, B^{(K)}, h^{(K)}, f^{(K)}, W^{(K)})$.
Then we have that\\
\noindent{\bf (i)} $W^{(K)}= \pi^{(K)}(W)$ is a controllability
submanifold of $SP$-reduced CH system $(M^{(K)}, B^{(K)}, h^{(K)},\\
f^{(K)}, W^{(K)})$ if and only if $W$ is $G$-invariant
controllability
submanifold of symmetric CH system $(T^\ast Q,G,B,H,F,W). $\\

\noindent{\bf (ii)} $D^{(K)}= T\pi^{(K)}(D^G) \subset
TM^{(K)}|_{W^{(K)}}$ is controllability distribution of the
$SP$-reduced CH system, such that $D_{W^{(K)}}=D^{(K)}\cap TW^{(K)}$
is a smooth regular integrable distribution on $W^{(K)}$, and the
presheaf $C_{W^{(K)}/D_{W^{(K)}}}^\infty$ has the
$(D^{(K)},D_{W^{(K)}})$-local extension property if and only if $D^G
\subset TT^*Q|_W$ is $G$-invariant controllability distribution of
the symmetric CH system, such that $D_W=D^G\cap TW$ is a smooth
regular $G$-invariant integrable distribution on $W$, and the
presheaf $C_{W/D_W}^\infty$ has the $(D^G,D_W)$-local extension
property.\\

\noindent{\bf (iii)} The singular Poisson reducible CH system
$(T^*Q, G, B, H, F, W)$ is Poisson reducible by $G$-invariant
controllability distribution $D^G$, if and only if the associated
$SP$-reduced CH system $(M^{(K)}, B^{(K)}, h^{(K)}, f^{(K)},
W^{(K)})$ is Poisson reducible by the reduced controllability
distribution $D^{(K)}$.
\end{theo}

\noindent{\bf Proof.} {\bf (i)} Because the symmetric CH system
$(T^\ast Q,G,B,H,F,W)$ is a singular Poisson reducible, for $u:
T^\ast Q\rightarrow W $, there exists $u^{(K)}: M^{(K)} \rightarrow
W^{(K)}$, such that $u^{(K)} \cdot \pi^{(K)}= \pi^{(K)}\cdot u $ and
\begin{equation*}X_{(M^{(K)}, B^{(K)}, h^{(K)},
f^{(K)}, u^{(K)})}\cdot \pi^{(K)}=T\pi^{(K)}\cdot X_{(T^\ast
Q,G,B,H,F,u)},
\end{equation*}
where $T\pi^{(K)}: TT^*Q \rightarrow TM^{(K)}$ is tangent map of the
projection $\pi^{(K)}: T^*Q \rightarrow M^{(K)}. $ Note that
$W^{(K)}= \pi^{(K)}(W)$, and $W$ is $G$-invariant,
hence from Definition 3.1 we know that {\bf (i)} holds.\\

\noindent{\bf (ii)} Notice that $D^{(K)}= T\pi^{(K)}(D^G)\subset
TM^{(K)}|_{W^{(K)}}, $ and $D^G \subset TT^*Q|_{W}$ is
$G$-invariant. Thus,
$$D_{W^{(K)}}=D^{(K)}\cap TW^{(K)} = T\pi^{(K)}(D^G) \cap T\pi^{(K)}(TW)
= T\pi^{(K)}(D^G \cap TW)= T\pi^{(K)}(D_{W}).$$ Since the bundle
structures of $T^*Q$ and $(T^\ast Q)/G$ is compatible with the
stratification of $T^*Q /G$, and $\pi^{(K)}$ and $T\pi^{(K)}$ are
vector bundle maps, and $\pi^{(K)}$ is Poisson, from Definition 2.4,
Definition 3.2 and the definition of the presheaf local extension
property, we know that {\bf (ii)} holds.\\

\noindent{\bf (iii)} If the singular Poisson reducible CH system
$(T^*Q, G, B, H, F, W)$ is Poisson reducible by $G$-invariant
controllability distribution $D^G$, we shall prove that the
associated $SP$-reduced CH system $(M^{(K)}, B^{(K)}, h^{(K)},
f^{(K)}, W^{(K)})$ is Poisson reducible by the reduced
controllability distribution $D^{(K)}$. By the above conclusions
{\bf (i), (ii)} and Theorem 3.4, it suffices to show that for any
$[z]= \pi^{(K)}(z) \in W^{(K)}, \; z\in W $ we have that
$(B^{(K)})^\sharp(\Delta^{(K)}_{[z]})\subset
(\Delta^{(K)}_{[z]})_{W^{(K)}}^\circ $, where .
\begin{equation*}\Delta^{(K)}_{[z]}:=\left\{\mathrm{d}f([z])\left|\begin{aligned}f\in
  C_{M^{(K)}}^\infty(U_{[z]}),\; \mathrm{d}f(y)|_{D^{(K)}(y)}=0,\; \mbox{for all } y\in
  U_{[z]}\cap W^{(K)},\\ \mbox{ and for any open neighborhood }U_{[z]} \mbox{ of } [z] \mbox{ in }
  M^{(K)} \end{aligned}\right.\right\},\end{equation*}
  and \begin{equation*}
    (\Delta^{(K)}_{[z]})_{W^{(K)}}:=\left\{\mathrm{d}f([z])\in \Delta^{(K)}_{[z]}\left|\begin{aligned}f|_{U_{[z]}\cap W^{(K)}}\
    \mbox{is constant for any open}\\ \mbox{ neighborhood } U_{[z]}
    \mbox{ of } [z] \mbox{ in } M^{(K)} \end{aligned}\right.\right\}.
\end{equation*}
In fact, assume that for any $\alpha^{(K)}=\mathrm{d}f_K([z])\in
\Delta^{(K)}_{[z]}, $ that is, $f_K \in C_{M^{(K)}}^\infty(U_{[z]}),
$ where $U_{[z]}$ is an open neighborhood of $[z]$ in $M^{(K)}, $
such that $U_{[z]}\cap W^{(K)}$ is $D_{W^{(K)}}$-invariant and $
\mathrm{d}f_K(y)|_{D^{(K)}(y)}=0,\; \mbox{for all } y\in U_{[z]}\cap
W^{(K)}. $ Since the map $\pi^{(K)}: T^*Q \rightarrow M^{(K)}$ is
surjective, then for $z \in T^*Q, $ such that $[z]= \pi^{(K)}(z), $
there is an open $G$-invariant neighborhood $\tilde{U}_z$ of $z$ in
$T^*Q, $ such that $U_{[z]}= \pi^{(K)}(\tilde{U}_{z})$ and
$\tilde{U}_{z}\cap W$ is $D_{W}$-invariant, and there exists a
$G$-invariant function $f \in C_{T^*Q}^\infty(\tilde{U}_{z})^G, $
such that $f= f_K\cdot (\pi^{(K)})$ and $
\mathrm{d}f(\tilde{y})|_{D^G(\tilde{y})}=0,\; \mbox{for all }
\tilde{y} \in \tilde{U}_{z}\cap W, $ and $\alpha=\mathrm{d}f(z)\in
\Delta^G_{z}. $ In the same way, for any
$\beta^{(K)}=\mathrm{d}g_K([z])\in (\Delta^{(K)}_{[z]})_{W^{(K)}}, $
that is, $g_K \in C_{M^{(K)}}^\infty(U_{[z]}) $ and
$g_K|_{U_{[z]}\cap W^{(K)}}$ is constant for any open neighborhood
$U_{[z]}$ of $[z]$ in $M^{(K)}$, then there is a $G$-invariant
function $g \in C_{T^*Q}^\infty(\tilde{U}_{z}), $ such that $g=
g_K\cdot (\pi^{(K)})$ and $g|_{\tilde{U}_{z}\cap W}$ is constant for
the corresponding open $G$-invariant neighborhood $\tilde{U}_{z}$ of
$z$ in $T^*Q$, and $\beta=\mathrm{d}g(z)\in (\Delta^G_{z})_{W}. $
Because the singular Poisson reducible CH system $(T^*Q, G, B, H, F,
W)$ is Poisson reducible by $G$-invariant controllability
distribution $D^G$, from the above conclusions {\bf (i), (ii)} and
Theorem 4.3 we have that $B^\sharp(\Delta^G_{z})\subset
(\Delta^G_{z})_{W}^\circ. $ It follows that
$$\{f, g\}_{B}(z)= <\mathrm{d}g(z),
B^\sharp(\mathrm{d}f(z))> = <\beta, B^\sharp(\alpha)> =0.$$ Notice
that the map $\pi^{(K)}: T^*Q \rightarrow M^{(K)}$ is Poisson, we
have that
$$0= \{f, g\}_{B}(z)= \{f_K\cdot (\pi^{(K)}), g_K\cdot
(\pi^{(K)})\}_{B}(z)= (\pi^{(K)})^{*} \{f_K,
g_K\}_{B^{(K)}}([z])=0.$$ Since $(\pi^{(K)})^{*}$ is injective, then
$\{f_K,g_K\}_{B^{(K)}}([z])=0, $ and $$ <\beta^{(K)},
(B^{(K)})^\sharp(\alpha^{(K)})> = <\mathrm{d}g_K([z]),
(B^{(K)})^\sharp(\mathrm{d}f_K([z]))>= \{f_K,g_K\}_{B^{(K)}}([z]) =
0. $$ Thus, we prove that $(B^{(K)})^\sharp(\alpha^{(K)})\in
(\Delta^{(K)}_{[z]})_{W^{(K)}}^\circ$ and
$(B^{(K)})^\sharp(\Delta^{(K)}_{[z]})\subset
(\Delta^{(K)}_{[z]})_{W^{(K)}}^\circ $, and hence the $SP$-reduced
CH system $(M^{(K)}, B^{(K)}, h^{(K)}, f^{(K)}, W^{(K)})$ is Poisson
reducible by the reduced controllability distribution $D^{(K)}$.\\

Conversely, if the associated $SP$-reduced CH system $(M^{(K)},
B^{(K)}, h^{(K)}, f^{(K)}, W^{(K)})$ is Poisson reducible by the
reduced controllability distribution $D^{(K)}$, by using the same
way we can verify that for any $z \in W, $ we have that
$B^\sharp(\Delta^G_{z})\subset (\Delta^G_{z})_{W}^\circ $. Thus, the
singular Poisson reducible CH system $(T^*Q, G, B, H, F, W)$ is
Poisson reducible by $G$-invariant controllability distribution
$D^G$ by the above conclusions {\bf (i), (ii)} and Theorem 4.3.
$\hskip 1cm \blacksquare$\\

Next, if considering the SPR-CH-equivalence of singular Poisson
reducible CH systems, we can get the following theorem to state the
property of Poisson reduction for singular Poisson reducible CH
systems by $G$-invariant controllability distribution leaves
invariant under SPR-CH-equivalence conditions.

\begin{theo}
Suppose that two singular Poisson reducible CH systems $(T^\ast
Q_i,G_i,B_i,H_i,F_i,W_i)$, $i=1,2,$ are SPR-CH-equivalent with
equivalent map $\varphi^*: T^\ast Q_2\rightarrow T^\ast Q_1$. Then
we have that\\

\noindent{\bf (i)} $W_1$ is $G_1$-invariant controllability
submanifold of CH system $(T^*Q_1, G_1, B_1, H_1, F_1, W_1)$ if and
only if $W_2$ is $G_2$-invariant controllability submanifold
of CH system $(T^*Q_2, G_2, B_2, H_2, F_2, W_2). $\\

\noindent{\bf (ii)} $D^{G_1}_1=T\varphi^*(D^{G_2}_2) \subset
TT^*Q_1|_{W_1}$ is $G_1$-invariant controllability distribution of
CH system $(T^*Q_1, G_1, B_1, H_1, F_1, W_1), $ such that
$D_{W_1}=D^{G_1}_1\cap TW_1$ is smooth regular $G_1$-invariant
integrable distributions on $W_1$ and the presheaf
$C_{W_1/D_{W_1}}^\infty$ has the $(D^{G_1}_1,D_{W_1})$-local
extension property, if and only if $D^{G_2}_2 \subset
TT^*Q_2|_{W_2}$ is $G_2$-invariant controllability distribution of
CH system $(T^*Q_2, G_2, B_2, H_2, F_2, W_2), $ such that
$D_{W_2}=D^{G_2}_2\cap TW_2$ is smooth regular $G_2$-invariant
integrable distributions on $W_2$ and the presheaf
$C_{W_2/D_{W_2}}^\infty$ has the $(D^{G_2}_2,D_{W_2})$-local
extension property.\\

\noindent{\bf (iii)} The singular Poisson reducible CH system
$(T^*Q_1, G_1, B_1, H_1, F_1, W_1)$ is Poisson reducible by
$G_1$-invariant controllability distribution $D^{G_1}_1$ if and only
if the singular Poisson reducible CH system $(T^*Q_2, G_2, B_2, H_2,
F_2, W_2)$ is Poisson reducible by $G_2$-invariant controllability
distribution $D^{G_2}_2$.
\end{theo}

\noindent{\bf Proof.} {\bf (i)} Because the singular Poisson
reducible CH systems $(T^\ast Q_i,G_i,B_i,H_i,F_i,W_i)$, $i=1,2,$
are SPR-CH-equivalent, then there exists a diffeomorphism $\varphi:
Q_1\rightarrow Q_2$, such that the cotangent lift map $\varphi^\ast=
T^\ast \varphi: T^\ast Q_2\rightarrow T^\ast Q_1$ and the tangent
lift map $T\varphi^\ast: TT^\ast Q_2 \rightarrow TT^\ast Q_1$ are
vector bundle isomorphism, and $\varphi^\ast$ is
$(G_2,G_1)$-equivariant Poisson map. Note that $W_1=\varphi^\ast
(W_2)$, and $W_i$ is $G_i$-invariant,$i=1,2$, hence from Definition
3.1 and Theorem 3.5 {\bf (i)} we know that {\bf (i)} holds.\\

\noindent{\bf (ii)} Notice that $D^{G_1}_1=T\varphi^*(D^{G_2}_2)
\subset TT^*Q_1|_{W_1}, $ and $D^{G_2}_2 \subset TT^*Q_2|_{W_2}$.
Thus,
$$D_{W_1}=D^{G_1}_1\cap TW_1 = T\varphi^*(D^{G_2}_2) \cap T\varphi^\ast
(TW_2)= T\varphi^\ast (D^{G_2}_2 \cap TW_2)= T\varphi^\ast
(D_{W_2}).$$ Since $\varphi^\ast$ and $T\varphi^\ast$ are vector
bundle isomorphism, and $\varphi^\ast$ is $(G_2,G_1)$-equivariant
Poisson map, from Definition 2.4, Definition 3.2, the definition of
the presheaf local extension
property and Theorem 3.5 {\bf (ii)}, we know that {\bf (ii)} holds.\\

\noindent{\bf (iii)} If singular Poisson reducible CH system
$(T^*Q_1, G_1, B_1, H_1, F_1, W_1)$ is Poisson reducible by
$G_1$-invariant controllability distribution $D^{G_1}_1$, from the
above conclusions {\bf (i), (ii)} and Theorem 6.1 we know that its
associated $SP$-reduced CH system $(M_{1}^{(K_1)}, B_{1}^{(K_1)},
h_{1}^{(K_1)}, f_{1}^{(K_1)}, W_{1}^{(K_1)})$ is Poisson reducible
by the reduced controllability distribution $D_{1}^{(K_1)}$, where
$D_{1}^{(K_1)}= T\pi_1^{(K_1)}(D_1^{G_1})$. On the other hand, from
Theorem 5.3 we know that two singular Poisson reducible CH systems
$(T^\ast Q_i,G_i,B_i,H_i,F_i,W_i)$, $i=1,2,$ are SPR-CH-equivalent
with equivalent map $\varphi^*: T^\ast Q_2\rightarrow T^\ast Q_1$,
if and only if the associated $SP$-reduced CH systems $(M_i^{(K_i)},
B_i^{(K_i)}, h_i^{(K_i)}, f_i^{(K_i)}, W_i^{(K_i)}), \; i=1,2,$ are
CH-equivalent with equivalent map $\varphi^*_{(K)}:
M_{2}^{(K_2)}\rightarrow M_{1}^{(K_1)}$. Moreover, from Theorem 3.5
we know that the property of Poisson reduction for CH systems by
controllability distribution leaves invariant under the
CH-equivalence. Thus, the $SP$-reduced CH system $(M_{2}^{(K_2)},
B_{2}^{(K_2)}, h_{2}^{(K_2)}, f_{2}^{(K_2)}, W_{2}^{(K_2)})$ is
Poisson reducible by the reduced controllability distribution
$D_{2}^{(K_2)}$, where $D_{1}^{(K_1)}
=T\varphi_{(K)}^\ast(D_{2}^{(K_2)}), $ and hence from the above
conclusions {\bf (i), (ii)} and Theorem 6.1 we have that the
singular Poisson reducible CH system $(T^*Q_2, G_2, B_2, H_2, F_2,
W_2)$ is Poisson reducible by $G_2$-invariant controllability
distribution $D^{G_2}_2$.\\

Conversely, if the singular Poisson reducible CH system $(T^*Q_2,
G_2, B_2, H_2, F_2, W_2)$ is Poisson reducible by $G_2$-invariant
controllability distribution $D^{G_2}_2$, by using the above same
argument we can get the desired conclusion. $\hskip 1cm
\blacksquare$\\

In particular, when $G$-action is free, the orbit space $T^*Q/G$ is
a smooth manifold, in this case we can consider the regular Poisson
reducible CH system and obtain the following corollary from the
above theorem and the regular Poisson reduction theorem given in
Wang and Zhang \cite{wazh12}.

\begin{coro}
Suppose that symmetric CH system $(T^\ast Q,G,B,H,F,W)$ is a regular
Poisson reducible CH system and its associated $RP$-reduced CH
system is $(T^\ast Q /G, B_{/G},h_{/G},f_{/G},W_{/G})$. Then we have
that\\
\noindent{\bf (i)} $W_{/G}= \pi_{/G}(W)$ is a controllability
submanifold of $RP$-reduced CH system $(T^\ast Q /G,
B_{/G},h_{/G},\\f_{/G},W_{/G})$ if and only if $W$ is $G$-invariant
controllability submanifold of
symmetric CH system $(T^\ast Q,G,B,H,F,W). $\\

\noindent{\bf (ii)} $D_{/G}= T\pi_{/G}(D^G) \subset T(T^\ast Q
/G)|_{W_{/G}}$ is controllability distribution of the $RP$-reduced
CH system, such that $D_{W_{/G}}=D_{/G}\cap TW_{/G}$ is a smooth
regular integrable distribution on $W_{/G}$, and the presheaf
$C_{W_{/G}/D_{W_{/G}}}^\infty$ has the $(D_{/G},D_{W_{/G}})$-local
extension property, if and only if $D^G \subset TT^*Q|_W$ is a
$G$-invariant controllability distribution, such that $D_W=D^G\cap
TW$ is a smooth regular $G$-invariant integrable distribution on
$W$, and the presheaf $C_{W/D_W}^\infty$ has the $(D^G,D_W)$-local
extension property.\\

\noindent{\bf (iii)} The regular Poisson reducible CH system $(T^*Q,
G, B, H, F, W)$ is Poisson reducible by $G$-invariant
controllability distribution $D^G$, if and only if the associated
$RP$-reduced CH system $(T^\ast Q /G, B_{/G},h_{/G},f_{/G},W_{/G})$
is Poisson
reducible by the reduced controllability distribution $D_{/G}$.\\
{\bf (iv)} If two regular Poisson reducible CH systems $(T^\ast
Q_i,G_i,B_i,H_i,F_i,W_i)$, $i=1,2,$ are RPR-CH-equivalent with
equivalent map $\varphi^*: T^\ast Q_2\rightarrow T^\ast Q_1$, then
CH system $(T^*Q_1, G_1, B_1, H_1, F_1, W_1)$ is Poisson reducible
by $G_1$-invariant controllability distribution $D^{G_1}_1$ if and
only if CH system $(T^*Q_2, G_2, B_2, H_2, F_2, W_2)$ is Poisson
reducible by $G_2$-invariant controllability distribution
$D^{G_2}_2$, where $D^{G_1}_1=T\varphi^*(D^{G_2}_2)$.
\end{coro}

\section{Applications}

In order to understand well the abstract theory, in this section,
some examples are given to state theoretical results of Poisson
reduction for CH systems by controllability distributions.\\

\noindent{\bf Example 7.1.} (Optimal point reduction of CH system)
Let $Q$ be a smooth manifold and $T^\ast Q$ its cotangent bundle
with the Poisson tensor $B$. Let $\Phi: G\times Q\rightarrow Q$ be a
smooth left action of the Lie group $G$ on $Q$, and the cotangent
lifted left action $\Phi^{T^\ast}: G\times T^\ast Q\rightarrow
T^\ast Q$ be canonical, and proper. Let $A_G'$ be the generalized
distribution on $T^\ast Q$ defined by the relation
$$A_G'(\alpha):=\{X_f(\alpha)|f\in C^\infty(T^\ast Q)^G\},\;\forall\,\alpha \in T^\ast Q, $$
and $A_G'$ be the $G$-characteristic distribution and it be smooth
and integrable in the sense of Stefan \cite{st74} and Sussmann
\cite{su73}. Moreover, let $\mathcal{J}: T^\ast Q \rightarrow T^\ast
Q/A'_G$ be the optimal momentum map associated to the cotangent
lifted left action $\Phi^{T^\ast}$, and $\rho\in T^\ast Q/A'_G$ be a
value of $\mathcal {J}$. If $m\in T^*Q $ and $\mathcal {J}(m)= \rho,
$ then $\mathcal {J}^{-1}(\rho) = A_G' \cdot m. $ Denote by $G_\rho$
the isotropy subgroup of $G$ in $\rho$. If $G_\rho$ acts properly on
$\mathcal {J}^{-1}(\rho)$, then the orbit space $(T^\ast
Q)_\rho=\mathcal{J}^{-1}(\rho)/G_\rho$ is a smooth symplectic
quotient manifold and the projection
$\pi_\rho:\mathcal{J}^{-1}(\rho)\rightarrow (T^\ast Q)_\rho$ is a
canonical surjective submersion. See Ortega and Ratiu \cite{orra04}.
For the symmetric CH system $(T^*Q, G, B, H, F, W)$, if the control
subset $W= \mathcal{J}^{-1}(\rho)$ is a $G$-invariant
controllability submanifold of $T^*Q$, and the generalized
distribution $D$ is given by the tangent spaces to the orbits of the
$G$-action. Assume that $D$ is $G$-invariant controllability
distribution of symmetric CH system. Notice that $g \cdot \mathcal
{J}(m)= \mathcal {J}(g \cdot m), $ for any $g \in G, \; m \in T^*Q,
$ this $G$-action is available on the leaf space of any distribution
spanned by $G$-equivariant vector fields. For any $\rho\in
T^*Q/A_G'$, there is a unique symplectic leaf $\mathcal{L}_\rho$ of
$(T^*Q, B)$ such that $\mathcal{J}^{-1}(\rho)\subset
\mathcal{L}_\rho$ and
$i_{\mathcal{L}_\rho}:\mathcal{J}^{-1}(\rho)\to \mathcal{L}_\rho$ is
smooth. In particular, since $(T^*Q, \omega)$ is a symplectic
manifold with a canonical symplectic form $\omega$ and the
$G$-action has an associated $\operatorname{Ad}^\ast$-equivariant
momentum map $\mathbf{J}: T^\ast Q \to \mathfrak{g}^\ast$, then the
fiber $W= \mathcal{J}^{-1}(\rho)$ of the optimal momentum map is
therefore an embedded submanifold of $T^*Q$. In this case, $D_W:=
D\cap TW $ is a generalized distribution given by the tangent spaces
to the orbits of the $G_\rho$-action, and $W/ D_W=
\mathcal{J}^{-1}(\rho)/G_\rho$, and from Ortega and Ratiu
\cite{orra04} we know that the presheaf $C_{W/D_W}^\infty$ has the
$(D,D_W)$-local extension property. Thus, the symmetric CH system
$(T^*Q, G, B, H, F, \mathcal{J}^{-1}(\rho))$ is Poisson reducible
by the controllability distribution $D$.\\

\noindent{\bf Example 7.2.} (Optimal orbit reduction of CH system)
Let $\rho\in T^\ast Q/A'_G$ be a value of the optimal momentum map
$\mathcal J$, and $\mathcal{O}_\rho=G\cdot \rho\subset T^\ast
Q/A'_G$ be the $G$-orbit of the point $\rho$. If isotropy subgroup
$G_\rho$ acts properly on $\mathcal J^{-1}(\rho)$, then there is a
unique smooth structure on $\mathcal J^{-1}(\mathcal{O}_\rho)$ that
makes it into an initial submanifold of $T^\ast Q$ and the
$G$-action on $\mathcal J^{-1}(\mathcal{O}_\rho)$ by restriction of
the $G$-action on $T^\ast Q$ is smooth and proper and all its
isotropy subgroups are conjugate to a given compact isotropy
subgroup of the $G$-action on $T^\ast Q$. Notice that
$G\cdot\mathcal J^{-1}(\rho)/G = \mathcal J^{-1}(\mathcal{O}_\rho)/G
$, the quotient space $(T^\ast Q)_{\mathcal{O}_\rho}=\mathcal
J^{-1}(\mathcal{O}_\rho)/G$ admits a unique smooth structure that
makes the projection $\pi_{\mathcal{O}_\rho}:\mathcal
J^{-1}(\mathcal{O}_\rho)\rightarrow (T^\ast Q)_{\mathcal{O}_\rho}$ a
surjective submersion. Moreover, note that the map
$\mathcal{J}^{-1}(\rho)/G_\rho\to
\mathcal{J}^{-1}(\mathcal{O}_\rho)/G$, $[m]_\rho\to
[m]_{\mathcal{O}_\rho}$, is a bijection, so the quotient space
$(T^*Q)_{\mathcal{O}_\rho}=\mathcal{J}^{-1}(\mathcal{O}_\rho)/G$ has
a smooth symplectic structure $\omega_{\mathcal{O}_\rho}$ induced
from the optimal point reduced space $(T^\ast Q)_{\rho}$. If the
$G$-action on $T^\ast Q$ is free, from the Poisson stratification
theorem given in Fernandes et al \cite{feorra09}, we know that the
symplectic leaves of the $RP$-reduced space $(T^\ast Q /G, B_{/G})$
are just given by the $O_O$-reduced space $((T^\ast
Q)_{\mathcal{O}_\rho},\omega_{\mathcal{O}_\rho})$, $ \rho\in T^\ast
Q/A'_G $. For the symmetric CH system $(T^*Q, G, B, H, F, W)$, if
the control subset $W= \mathcal{J}^{-1}(\mathcal{O}_\rho)$ is a
$G$-invariant controllability submanifold of $T^*Q$, and the
generalized distribution $D$ is given by the tangent spaces to the
orbits of the $G$-action. Assume that $D$ is $G$-invariant
controllability distribution of symmetric CH system. For $\rho\in
T^*Q/A_G'$, assume that $\mathcal{J}^{-1}(\mathcal{O}_\rho)$ is
connected and there is a unique symplectic leaf
$\mathcal{L}_{\mathcal{O}_\rho}$ of $(T^*Q, B)$ such that
$\mathcal{J}^{-1}(\mathcal{O}_\rho)\subset
\mathcal{L}_{\mathcal{O}_\rho}$ and
$i_{\mathcal{L}_{\mathcal{O}_\rho}}:\mathcal{J}^{-1}(\mathcal{O}_\rho)\to
\mathcal{L}_{\mathcal{O}_\rho}$ is smooth. In particular, since
$(T^*Q, \omega)$ is a symplectic manifold with a canonical
symplectic form $\omega$ and the $G$-action has an associated
$\operatorname{Ad}^\ast$-equivariant momentum map $\mathbf{J}:
T^\ast Q \to \mathfrak{g}^\ast$, then the fiber submanifold $W=
\mathcal{J}^{-1}(\mathcal{O}_\rho)$ of the optimal momentum map is
therefore an embedded submanifold of $T^*Q$. In this case, $D_W:=
D\cap TW $ is a generalized distribution given by the tangent spaces
to the orbits of the $G$-action on $W$, and $W/ D_W=
\mathcal{J}^{-1}(\mathcal{O}_\rho)/G $, and from Ortega and Ratiu
\cite{orra04} we know that the presheaf $C_{W/D_W}^\infty$ has the
$(D,D_W)$-local extension property. Thus, the symmetric CH system
$(T^*Q, G, B, H, F, \mathcal{J}^{-1}(\mathcal{O}_\rho))$ is Poisson
reducible by the controllability
distribution $D$.\\

\noindent{\bf Example 7.3.} (Poisson reduction of CH system by
characteristic distribution) For the CH system $(T^*Q, B, H, F, W)$,
assume that the control subset $W$ is an embedded controllability
submanifold of $T^*Q$, and the distribution $D:=
B^\sharp((TW)^\circ)\subset TT^*Q |_W $ is controllability
distribution, such that $D_W:= D\cap TW $ is a smooth and integrable
generalized distribution on $W$, and the presheaf $C_{W/D_W}^\infty$
has the $(D,D_W)$-local extension property. Then from Theorem 2.8,
we know that $(T^*Q,B,D,W)$ is Poisson reducible, and hence the CH
system $(T^*Q, B, H, F, W)$ is Poisson reducible by the
characteristic distribution $D:=
B^\sharp((TW)^\circ)$.\\

\noindent{\bf Example 7.4.} (Poisson reduction of CH system by
controllability coisotropic submanifold) For the CH system $(T^*Q,
B, H, F, W)$, assume that the control subset $W$ is an embedded
controllability coisotropic submanifold of $T^*Q$, and $D:=
B^\sharp((TW)^\circ)\subset TT^*Q |_W $ is controllability
distribution, such that $D_W:= D\cap TW $ and the presheaf
$C_{W/D_W}^\infty$ has the $(D,D_W)$-local extension property. In
this case, by using the coisotropic property of $W$ we can prove
that $B^\sharp((TW)^\circ)\subset TW, $ and $D_W:= D\cap TW =D $ is
a smooth and integrable generalized distribution on $W$. Then from
Theorem 2.8, we know that $(T^*Q,B,D,W)$ is Poisson reducible, and
hence the CH system $(T^*Q, B, H, F, W)$ is Poisson reducible by the
controllability coisotropic submanifold.\\

\noindent{\bf Example 7.5.} (Poisson reduction of CH system by the
inverse image of a coadjoint orbit) Let $(T^*Q, B)$ be a Poisson
fiber bundle with associated Poisson tensor $B \in
\wedge^2(T^*T^*Q)$. Let $G$ be a Lie group acting freely and
canonically on $(T^*Q, B)$ with an associated coadjoint equivariant
standard momentum map $\mathbf{J}: T^*Q \rightarrow \mathfrak{g}^*.
$ Let $\mu \in \mathfrak{g}^*$ be a value of the momentum map
$\mathbf{J}$ and $\mathcal{O}_\mu \subset \mathfrak{g}^*$ its
coadjoint orbit. We know that $\mathbf{J}^{-1}(\mathcal{O}_\mu)$ is
an initial submanifold of $T^*Q$, and for any $m\in
\mathbf{J}^{-1}(\mathcal{O}_\mu),$
$$T_m(\mathbf{J}^{-1}(\mathcal{O}_\mu))= ker T_m \mathbf{J} + \mathfrak{g}\cdot m. $$
If $\mathcal{L}$ is the symplectic leaf of $(T^*Q, B)$ containing
$m$, from Proposition 2.1 and $\mathfrak{g}\cdot m \subset T_m
\mathcal{L}$, we have that
$$B^\sharp((T_m(\mathbf{J}^{-1}(\mathcal{O}_\mu)))^\circ)\subset T_m (\mathbf{J}^{-1}(\mathcal{O}_\mu)).$$
Thus, $\mathbf{J}^{-1}(\mathcal{O}_\mu)$ is a coisotropic
submanifold of $(T^*Q, B)$. For the CH system $(T^*Q, B, H, F, W)$,
if the control subset $W= \mathbf{J}^{-1}(\mathcal{O}_\mu)$ is an
embedded controllability coisotropic submanifold of $T^*Q$, and $D:=
B^\sharp((TW)^\circ)\subset TT^*Q |_W $ is controllability
distribution, then $D_W:= D\cap TW =D $ is a smooth and integrable
generalized distribution on $W$. If the presheaf $C_{W/D_W}^\infty$
has the $(D,D_W)$-local extension property, then the CH system
$(T^*Q, B, H, F,\mathbf{J}^{-1}(\mathcal{O}_\mu))$ is Poisson
reducible by the inverse image of a
coadjoint orbit.\\

\noindent{\bf Example 7.6.} (Poisson reduction of CH system by
controllability cosymplectic submanifold) For the CH system $(T^*Q,
B, H, F, W)$, assume that the control subset $W$ is an embedded
controllability cosymplectic submanifold of $T^*Q$, and $D:=
B^\sharp((TW)^\circ)\subset TT^*Q |_W $ is controllability
distribution. In this case, by using the cosymplectic property of
$W$ we have that $B^\sharp((TW)^\circ) \cap TW = \{0\}, $ and hence
$D_W:= D\cap TW = \{0\}$ is a trivially smooth and integrable
generalized distribution on $W$, and the presheaf $C_{W/D_W}^\infty=
C_W^\infty $ has the $(D,D_W)$-local extension property. Then from
Theorem 2.8, we know that $(T^*Q,B,D,W)$ is Poisson reducible, and
hence the CH system $(T^*Q, B, H, F, W)$ is Poisson reducible by the
controllability cosymplectic submanifold.\\

The theory of mechanical control system is a very important subject,
following the theoretical development of geometric mechanics, a lot
of important problems about this subject are being explored and
studied. Wang in \cite{wa13d} studies the Hamilton-Jacobi theory of
regular controlled Hamiltonian systems with the symplectic structure
and symmetry, and Wang in \cite{wa13f, wa13e} apply the above work
to give explicitly the motion equations and Hamilton-Jacobi
equations of reduced rigid spacecraft-rotor system and reduced
underwater vehicle-rotors system on the symplectic leaves by
calculation in detail, which show the effect on controls in regular
symplectic reduction (by stages) and Hamilton-Jacobi theory. But if
we define a controlled Hamiltonian system on the cotangent bundle
$T^*Q$ by using a Poisson structure, just same as we have done in
this paper and in Wang and Zhang \cite{wazh12}, then the way given
in Wang \cite{wa13d} cannot be used, what and how we could do? This
is a problem worthy to be considered in detail. In addition, we also
note that there have been a lot of beautiful results of reduction
theory of Hamiltonian systems in celestial mechanics, hydrodynamics
and plasma physics. Thus, it is an important topic to study the
application of reduction theory and Hamilton-Jacobi theory of
controlled Hamiltonian systems in celestial mechanics, hydrodynamics
and plasma physics. These are our goals in
future research.\\

\noindent{\bf Acknowledgments:} The authors would like to thank the
support of Nankai University, 985 Project and the Key Laboratory of
Pure Mathematics and Combinatorics, Ministry of Education, China.\\


\begin{thebibliography}{99}

\bibitem{abma78}
R. Abraham and J. E. Marsden, Foundations of Mechanics, Second
edition, Addison-Wesley, Reading, MA, 1978.
\bibitem{abmara88}
R. Abraham, J. E. Marsden and T. S. Ratiu, Manifolds, Tensor
Analysis and Applications, Applied Mathematical Science, 75,
Springer-Verlag, New York, 1988.
\bibitem{bipura04}
P. Birtea, M. Puta and T. S. Ratiu, Controllability of Poisson
systems, SIAM J. Control Optim., 43(3) (2004), 937-954.
\bibitem{blle02}
A.M. Bloch and N.E. Leonard, Symmetries, conservation laws, and
control, In ``Geometry, Mechanics and Dynamics, Volume in Honor of
the 60th Birthday of J.E. Marsden" (eds. P.Newton, P.Holmes and A.
Weinstein), Springer, New York, 2002.
\bibitem{faza08}
F. Falceto and M. Zambon, An extension of the Marsden-Ratiu
reduction for Poisson manifolds, Lett. Math. Phys., 85(3)(2008),
203-219.
\bibitem{feorra09}
R. L. Fernandes, J. P. Ortega and T. S. Ratiu, The momentum map in
Poisson geometry, Amer. J. Math., 131(5) (2009), 1261-1310.
\bibitem{jora09}
M. Jotz and T. S. Ratiu, Poisson reduction by distributions, Lett.
Math. Phys., 87(1-2) (2009), 139-147.
\bibitem{jorasn11}
M. Jotz, T. S. Ratiu, and J. $\acute{S}$niatycki, Singular reduction
of Dirac structures, Trans. Amer. Math. Soc., 363 (2011), 2967-3013.
\bibitem{lema97}
N.E. Leonard, and J.E. Marsden, Stability and drift of underwater
vehicle dynamics: mechanical systems with rigid motion symmetry,
Physica D 105(1997), 130-162.
\bibitem{lima87}
P. Libermann and C. M. Marle, Symplectic Geometry and Analytical
Mechanics, Kluwer Academic Publishers, 1987.
\bibitem{ma92}
J. E. Marsden, Lectures on Mechanics, in: London Mathematical
Society Lecture Notes Series, vol. 174, Cambridge University Press,
1992.
\bibitem{mamiorpera07}
J. E. Marsden, G. Misiolek, J. P. Ortega, M. Perlmutter and T. S.
Ratiu, Hamiltonian Reduction by Stages, in: Lecture Notes in
Mathematics, vol. 1913, Springer, 2007.
\bibitem{mara86}
J. E. Marsden and T. S. Ratiu, Reduction of Poisson manifolds, Lett.
Math. Phys., 11(2) (1986), 161-169.
\bibitem{mara99}
J. E. Marsden and T. S. Ratiu, Introduction to Mechanics and
Symmetry, second edition, Texts in Applied Mathematics, vol. 17,
Springer-Verlag, 1999.
\bibitem{mawazh10}
J. E. Marsden, H. Wang, and Z. X. Zhang, Regular reduction of
controlled Hamiltonian system with symplectic structure and
symmetry, (arXiv: 1202.3564, To appear in Diff. Geom. Appl.).
\bibitem{mawe74}
J. E. Marsden and A. Weinstein, Reduction of symplectic manifolds
with symmetry, Rep. Math. Phys., 5 (1974), 121-130.
\bibitem{nivds90}
H. Nijmeijer and A. J. Van der Schaft, Nonlinear Dynamical Control
Systems, Springer-Verlag, 1990.
\bibitem{orra04}
J. P. Ortega and T. S. Ratiu, Momentum Maps and Hamiltonian
Reduction, Progress in Mathematics, 222, Birkh\"{a}user, 2004.
\bibitem{pf01}
M. J. Pflaum, Analytic and Geometric Study of Stratified Spaces,
Lecture Notes in Mathematics, 1768, Springer-Verlag, 2001.
\bibitem{sjle91}
R. Sjamaar and E. Lerman, Stratified symplectic spaces and
reduction, Ann. of Math., 134 (1991), 375-422.
\bibitem{sn06}
J. $\acute{S}$niatycki, Singular reduction for nonlinear control
systems, Reports on Methematical Physics, 57(2) (2006), 163-178.
\bibitem{sa89}
G. S$\acute{a}$nchez de Alvarez, Controllability of Poisson control
systems with symmetry, Contemp. Math., 97 (1989), 399-412.
\bibitem{st74}
P. Stefan, Accessible sets, orbits and foliations with
singularities, Proc. Lond. Math. Soc. 29 (1974), 699-713.
\bibitem{su73}
H. Sussmann, Orbits of families of vector fields and integrability
of distribution, Trans. Amer. Math. Soc. 180 (1973), 171-188.
\bibitem{wa12}
H. Wang, Singular reduction of regular controlled Hamiltonian system
with symmetry, (2012).
\bibitem{wa13d}
H. Wang, Hamilton-Jacobi theorems for regular controlled Hamiltonian
system and its reductions, (2013d, arXiv: 1305.3457, To submit to J.
Geom. Mech. ).
\bibitem{wa13f}
H. Wang, Symmetric reduction and Hamilton-Jacobi equation of rigid
spacecraft with a rotor, J. Geom. Symm. Phys., 32 (2013), 87-111,
(arXiv: 1307.1606 ).
\bibitem{wa13e}
H. Wang, Symmetric reduction and Hamilton-Jacobi equation of
underwater vehicle with internal rotors, (2013e, arXiv: 1310.3014 ).
\bibitem{wazh12}
H. Wang and Z. X. Zhang, Optimal reduction of controlled Hamiltonian
system with Poisson structure and symmetry, Jour. Geom. Phys., 62
(5) (2012), 953-975.
\bibitem{we83}
A. Weinstein, The local structure of Poisson manifolds, Jour. Diff.
Geom., 18 (1983), 523-557.

\end{thebibliography}
\end{document}